\documentclass[12pt]{amsart}
\usepackage{amssymb,amsmath,amscd}
\usepackage{latexsym,bm,bbm,mathrsfs}
\usepackage{hyperref,graphicx}
\usepackage{mathtools}
\usepackage{stmaryrd}
\usepackage[all,pdf]{xy}
\usepackage{graphicx}
\usepackage{tikz}
\usepackage{tikz-cd}
\usepackage{fontawesome5}
\tikzset{
  symbol/.style={
    draw=none,
    every to/.append style={
      edge node={node [sloped, allow upside down, auto=false]{$#1$}}}
  }
}
\newcommand*\tcircle[1]{%
  \raisebox{-0.5pt}{%
    \textcircled{\fontsize{8.5pt}{0}\fontfamily{phv}\selectfont #1}%
  }%
}
\graphicspath{ {./} }

\newtheorem{thm}{Theorem}[section]

\newtheorem{example}[thm]{Example}
\newtheorem{theorem}[thm]{Theorem}
\newtheorem{cor}[thm]{Corollary}
\newtheorem{corollary}[thm]{Corollary}
\newtheorem{prop}[thm]{Proposition}
\newtheorem{lemma}[thm]{Lemma}

\newtheorem{remark}[thm]{Remark}

\newtheorem{defn}[thm]{Definition}

\newcommand{\bbp}{{\mathbb{P}}}

\newcommand{\bbz}{{\mathbb{Z}}}

\newcommand{\bbf}{\mathbb{F}}

\newcommand{\Tr}{\operatorname{Tr}}

\newcommand{\End}{\operatorname{End}}

\title{Waring Problem for Matrices over Finite Fields}

\author{Krishna Kishore, Adrian Vasiu, and Sailun Zhan}

\thanks{\faCreativeCommons\faCreativeCommonsBy\faCreativeCommonsNc\faCreativeCommonsNd This manuscript is made available under the CC-BY-NC-ND 4.0 license.}
\thanks{https://creativecommons.org/licenses/by-nc-nd/4.0/ }
\thanks{https://doi.org/10.1016/j.jpaa.2024.107656}


\thanks{The first author is supported by Science and Engineering Research Board (SERB) MATRICS grant MTR/2021/000319 of the Government of India.
}

\begin{document}
\begin{abstract}
 We prove that for all integers $k \geq 1$, $q\ge (k-1)^4+ 6k$, and $m \geq 1$, every matrix in $ M_m(\mathbb F_q)$ is a sum of two kth powers: $M_m(\mathbb F_q)=\{A^k+B^k|A,B\in M_m(\mathbb F_q)\}$. We further generalize and refine this result in the cases when both $B$ and $C$ can be chosen to be invertible, cyclic, or split semisimple, when $k$ is coprime to $p$, or when $m$ is sufficiently large. We also give a criterion for the Waring problem in terms of stabilizers.

\end{abstract}

\maketitle

\section{Introduction}
The Waring problem for matrices is to address whether matrices over a ring can be expressed as a sum of two $k$th powers of matrices. We begin by establishing a framework in which the {\it Waring matrix problem} can be formulated precisely. Let $l,k,m$ and $s$ be positive integers. Let $\mathcal S$ be a semiring with identity and  $\mathcal T$ be a subset of $\mathcal S$. Let $$s\mathcal T:=\{\sum_{i=1}^s x_i|x_1,\ldots,x_s\in\mathcal T\},$$the sum of $s$ elements in $\mathcal T$; in particular $1\mathcal T=\mathcal T$. Let $M_m(\mathcal S)$ be the semiring of $m\times m$ matrices with entries in $\mathcal S$.

Consider the set
$$\mathbb P_{k,m,\mathcal S}:=\{A^k|A\in M_m(\mathcal S)\},
$$
of $k$th powers of matrices in $M_m(\mathcal S)$.

In analogy with the case of the semiring $Z_{\ge 0}$ of nonnegative integers over which the classical Waring problem is formulated, the Waring matrix problem for a class $\mathcal C$ of commutative semirings would, in the most refined form, asks to classify all quadruples $(m,\mathcal S,k,s)$, with $\mathcal S$ in the class $\mathcal C$, such that 
$$M_m(\mathcal S)=s\mathbb P_{k,m,\mathcal S};$$ 
if $s$ is minimal with this property, then it is called the $k$th Pythagoras number of $M_m(\mathcal S)$ in \cite[Section 1.2]{Pum07}. One can formulate different generalizations of the problem by considering only $k$th powers of invertible, idempotent,  (split) semisimple, or cyclic matrices.\footnote{ Recall that a matrix $A$ in $M_n(K)$, where $K$ is a field, is called {\it semisimple} if the minimal polynomial of $A$ is a product of distinct monic irreducible and separable polynomials; if moreover all these irreducible polynomials have degree $1$, then $A$ is called {\it split semisimple} or {\it diagonalizable}. }

In this paper we study the class of finite fields over which there exist extensive results on the number of different concrete types of square matrices; see the survey \cite{Mor06}. Several partial results to the Waring matrix problem over finite fields have been obtained recently; see \cite[Theorem 1.1]{TZS19}, \cite[Theorems 4.2 and 4.3]{DK19}, \cite{Kis22}, and \cite{KS23}. One of our results gives an alternative proof for the one proven by Kishore and Singh \cite{KS23} with an explicit lower bound on $q$. In other words, while Kishore and Singh established that for every $k \geq 1$ there exists a constant $C(k)$ depending only on $k$ such that for all $q \ge C(k)$ for all $n \geq 1$ every matrix in $M_n(\mathbb F_q)$ is a sum of two $k$th powers, the lower bound $C(k)$ was not given explicitly. In this article, we give an explicit lower bound. To state these partial results as well as our new results, we introduce some notation. Let $p$ be a prime, and  $q:=p^l$ and let $n\in\mathbb N\setminus\{1\}$. Define the following sets:
\[
M_{n,q}:=M_n(\bbf_q),\ G_{n,q}:=GL_n(\bbf_q),\ S_{n,q}:=SL_n(\bbf_q),
\] 
\[Z_{n,q}:=\{aI_n|a\in\mathbb F_q,\ a^n=1\}\ {\rm is\ the\ center\ of\ }S_{n,q},
\]
\[
s\mathbb P_{k,m,q}:=s\mathbb P_{k,m,\mathbb F_q}, \mathbb P_{k,q}:=\mathbb P_{k,1,q}\subset\mathbb F_q,
\]
\[P=P_{k,n,q}:=2\mathbb P_{k,n,q}=\{A^k+B^k\mid A,B\in M_{n,q}\},
\]
\[
Q=Q_{k,n,q}:=\{A^k+B^k\mid A,B\in G_{n,q} \},
\]
\[
P^{ss}=P^{ss}_{k,n,q}:=\{A^k+B^k\mid {\rm split\ semisimple}\ A,B\in M_{n,q} \},
\]
\[
Q^{ss}=Q^{ss}_{k,n,q}:=\{A^k+B^k\mid {\rm split\ semisimple}\ A,B\in G_{n,q} \},
\]
\[\mathbb I_{n,q}:=\{A\in M_{n,q}|A^2=A\},\ {\rm and}\ \Pi_{n,q}:=2\mathbb I_{n,q}.
\]

We will study the following inclusions
\begin{center}
\begin{tikzcd}[column sep=.7em]
&Q_{k,n,q}^{ss}\arrow[d,symbol=\subset] \arrow[r,symbol=\subset] &Q_{k,n,q} \arrow[d,symbol=\subset] \\
\Pi_{n,q}\arrow[r,symbol=\subset]& P_{k,n,q}^{ss} \arrow[r,symbol=\subset] & P_{k,n,q} \arrow[r,symbol=\subset] & M_{n,q},
\end{tikzcd}
\end{center}
such as when some of them are equalities. E.g., from \cite[Lemma 3]{HT88} it follows that always $\Pi_{n,q}\subsetneq M_{n,q}$. If $q=2$, then for $M_{n,2}$ the notions idempotent and split semisimple coincide, hence $P_{k,n,2}^{ss}=\Pi_{n,2}\subsetneq M_{n,2}$.

With the above notations, \cite[Theorem 1.1]{TZS19} shows that for $q\in\{2,3\}$, each matrix in $M_{n,q}$ is a sum of three idempotents and hence every $k$th
Pythagoras number of $M_{n,q}$ is either $2$ or $3$ and the remaining problem is to determine all pairs $(n,k)$ when it is exactly 2 (i.e., when $P=M_{n,q}$). The result \cite[Theorem 4.2]{DK19} implies that if $k\in\{3,5,7,9,11\}$ and $p$ does not belong to a finite set $E_k$ of primes, then $P_{k,n,q}=M_{n,q}$ (e.g., $E_3=\{2,3,7\}$ and $E_5=\{2,5,11,31,41,61\}$). Also, \cite[Theorem 1.1]{KS23} shows the existence of a constant $C(k)\in\mathbb N$ such that for $q\ge C(k)$ we have $P=M_{n,q}$, but no explicit value of $C(k)$ was provided and the methods used in loc. cit. do require $C(k)$ to be much bigger than the bounds one expects to get using the results on when $2\mathbb P_{k,q}$ is $\mathbb F_q$ obtained in \cite{Sma77} based on \cite{Jol73}. 

We could not find results in the literature pertaining to $Q$, $P^{ss}$, or $Q^{ss}$ besides the identity $M_{n,q}=Q_{1,n,q}$ proved in \cite[Theorem 12 or Corollary 13]{Hen74} or \cite{Lor87}, the identity $Q_{2,n,q}=M_{n,q}$ (resp. inclusion $Q_{2,n,q}\supset M_{n,q}\setminus\{I_n\}$) for $q>2$ and $n$ even (resp. and $n$ odd) proved in \cite[Theorem 1]{Yo98}, the identity $P^{ss}_{1,n,q}= M_{n,q}$ for $q\ge 3$ proved in \cite[Theorem 2.6]{Bot00} and \cite[Lemma 2]{NVZ04} which easily implies that $P_{k,n,q}=M_{n,q}$ if $\gcd(k,w_{n,q})=1$, where $w_{n,q}$ is the product of primes that divide the exponent $e_{n,q}$ of $G_{n,q}$.

If the pair $(n,q)$ is fixed, there exist infinitely many $k$s such that $P_{k,n,q}\subsetneq M_{n,q}$. E.g., if $e_{n,q}$ divides $k$, then $\mathbb P_{k,n,q}=\mathbb I_{n,q}$ (see Proposition \ref{exponent}), hence $P_{k,n,q}=\Pi_{n,q}\subsetneq M_{n,q}$. If $q+n\le 5$, in Sections \ref{S6} and \ref{S7} we check that we have a strict inclusion $P_{k,n,q}\subsetneq M_{n,q}$ if and only if $w_{n,k}$ divides $k$. Like $w_{2,2}=e_{2,2}=6$, $w_{3,2}=42$ and $e_{3,2}=84$, and $w_{2,3}=6$ and $e_{2,3}=24$ (see Theorems \ref{M22}, \ref{M32}, and \ref{M23}). 

As it is not clear to us how to extend the last results  to all possible pairs $(n,q)$, in what follows we will mainly concentrate on getting results that either are true for all $n\ge 2$ or require some of the following types of conditions to hold: (i) $q$ is greater than expressions in suitable divisors of $k$, (ii) $n\gg 0$, or (iii) certain coprime conditions hold.

Define $d_m:=\gcd(k,q^m-1)$. We will abbreviate $d:=d_1$. For all $m\in\mathbb N$ we have an identity 
\begin{equation}\label{EQ1}
\mathbb P_{k,q^m}=\mathbb P_{d_m,q^m}.
\end{equation}
For all pairs $(k,n)$ as above there exists a smallest natural number $\gamma(k,q)=\gamma(d,q)$ and a divisor $\ell$ of $l$ such that we have

\begin{equation}\label{EQ2}
\bigcup_{s=1}^{\infty} s\mathbb{P}_{k,q}=\gamma(d,q)\mathbb{P}_{d,q}=d\mathbb{P}_{d,q}=\bbf_{p^{\ell}}
\end{equation}
(see \cite[Lemma 1 and Theorem 1]{Tor38}), so $1\le\gamma(d,q)\le d$, and these results have a great impact on the case $n\ge 2$ as well.\footnote{If $\ell<l$, our notation $\gamma(d,q)$ differs from the common notation in the literature which simply defines $\gamma(d,q)$ to be $\infty$ and hence which loses some information.}

In this paper we study different Waring problems for many nonempty subsets $S\subset M_{n,q}$ that are of interest, such as $P$ or $Q$. For each such $S$ we define two subgroups of $G_{n,q}$ as follows:
\[
L_S:=\{g\in G_{n,q}\mid gS=S\}\ {\rm and}\ R_S:=\{g\in G_{n,q}\mid Sg=S\}.
\]

Different properties of $S$, such as being stable under conjugation or transposition or scalar multiplication or addition, are closely related to properties of $L_S$ and/or $R_S$. For instance, if $d=1$, then $P$ and $P^{ss}$ are stable under scalar multiplication, Similarly, if $S$ is stable under conjugation and an $\mathbb F_q$-vector subspace of $M_{n,q}$, then for $(n,q)\notin\{(2,2),(2,3)\}$, from the representation results \cite[Propositions 3.3 and 3.4 and Theorem 3.10 1) and 2)]{Vas03} it follows that there exist only 4 possibilities for $S$: either $S\in\{M_{n,q},0\}$ and $L_S=R_S=G_{n,q}$ or $S\in\{\mathbb F_qI_n,V_{n,q}\}$ and $L_S=R_S$ is the center of $G_{n,q}$, where $V_{n,q}:=\{A\in M_{n,q}|\Tr(A)=0\}$. In practice $S$ is stable under conjugation and different from $0$, $\mathbb F_qI_n$, and $V_{n,q}$, but it is hard to decide if $S$ is or is not an $\mathbb F_q$-vector subspace. But for $P$ and $Q$, in Section \ref{S2} we prove the following simpler criterion:

\begin{theorem}\label{thm1}
We assume that $(n,q)$ is neither $(2,2)$ nor $(2,3)$. Let $g\in S_{n,q}\setminus Z_{n,q}$. Then the following properties are true:

(a) The equality $P=M_{n,q}$ holds if and only if $g\in L_P$.

(b) If $p\nmid k$ and $-1\in\mathbb P_{d,q}$, then $Q=M_{n,q}$ if and only if $g\in L_Q$. 

(c) If $p\nmid k$, then $Q\cup\{0_n\}=M_{n,q}$ if and only if $g\in L_Q$. 
\end{theorem}

Using lower-upper triangular decompositions $A=L+U$ and block decompositions, in Section \ref{S4} we prove the following theorems:

\begin{theorem}\label{thm2} 
Suppose $d=1$. Then the following properties hold:

(a) If $q=3$ and $3$ does not divide $k$, then $Q=M_{n,q}$.

(b) If $q\ge 3$, then $P^{ss}=M_{n,q}$. 

(c) If $q\ge 4$, then $Q^{ss}=M_{n,q}$.
\end{theorem}

Theorem \ref{thm2} (b) for $k=1$ was first proved in \cite[Theorem 2.6]{Bot00}. When $q=2$, we also prove the case for $P$ (resp. $Q$) when $k$ is odd (resp. $\gcd(k,6)=1$); see Theorems \ref{thm3} and \ref{thm4}.

\begin{defn}
For a commutative ring $\mathcal R$, we say that $A\in M_{m}(\mathcal R)$ is \emph{cyclic} if, under the left multiplication action of $A$ on $\mathcal R^m$, $\mathcal R^m$ is a cyclic $\mathcal R[A]$-module, i.e., there exists an element $v\in\mathcal R^m$ such that $\{v,Av,\ldots,A^{m-1}v\}$ generates $\mathcal R^m$ (equivalently, is an $\mathcal R$-basis of $\mathcal R^m$).
\end{defn}

Let $K$ be a field. If $A\in M_m(K)$, then $A$ is cyclic if and only if $\mu_A(x)=\chi_A(x)$.

Using a counting argument (resp. lower-upper triangular decompositions and block decompositions) in Theorem \ref{thm5} we prove that if $q\ge3$, then each matrix in $M_{n,q}$ is a sum of two invertible cyclic (resp. invertible semisimple) matrices.

If $n=2$ and $P^{ss}=M_{2,q}$, then by reasons of traces each element of $\mathbb F_q$ is a sum of four $k$th powers and this implies that the inequality $\sum_{i=1}^4 (q-1)^id^{-i}\ge q-1$ holds. Thus if
$q<1+\sqrt[3]{\frac{1}{4}}d^{\frac{4}{3}}$, then for $n=2$ we have $P^{ss}\neq M_{2,q}$. Moreover, the current results on the Waring problem for finite fields guarantee that each element of $\mathbb F_q$ is a sum
of four $k$th powers only if $q>(d-1)^{\frac{8}{3}}$ (see Theorem \ref{Joly} for this result due to Joly). In order to get a converse for $n=2$ and to move from $n=2$ to all integers $n\ge 2$, our methods require increasing the exponent $\frac{8}{3}$ to $4$ as follows.

\begin{theorem}\label{thm7}
If $q\ge (d-1)^4+6d$, then $P^{ss}=M_{n,q}$ for all $n\ge 2$.
\end{theorem}

As $(d-1)^4+6d\le (k-1)^4+6k$, the above theorem proved in Section \ref{S4} gives a new proof of \cite[Theorem 1.1]{KS23} in a stronger form that involves $P^{ss}$ and not $P$ and with the explicit bound $C(k)=(k-1)^4+6k$.

We study the relation between nonscalar matrices and quasi-cyclic matrices (see Definition \ref{defqcyc}) in Theorem \ref{nonscalar}. Using this as a main new tool and results on the Waring problem over finite fields, i.e., Equation (\ref{EQ2}) and different refinements of the inequality $\gamma(d,q)\le d$ (see \cite{Jol73}, \cite{Sma77}, \cite{Win98}, \cite{Cip09}, \cite{CC12}, etc.),  in Section \ref{S4} we will prove the following series of results.

\begin{theorem}\label{thm8}
Suppose $p$ does not divide $k$. Then $P=M_{n,q}$ if 
\[
q\ge \max\{\sqrt{(d-1)^6+6d},(d_2-1)^2+1,\lfloor (d_3-1)^{4/3}\rfloor+1\}.
\]
\end{theorem}

As $d_m\le k$ for all $m\in\mathbb N$, from Theorem \ref{thm8} we get directly: 

\begin{cor}\label{cor8}
If $p\nmid k$ and $q\ge k^3-3k^2+3k$, then $P=M_{n,q}$.
\end{cor}

Corollary \ref{cor8} strengthens \cite[Theorem 4.3]{DK19} for the case of the finite rings $M_{n,q}$ as loc. cit. required $k$ to be odd and $q>(k-1)^4$. 

In an email to the first author dated July 7, 2021, Michael Larsen conjectured that for a given $k$, there exists a constant $N(k)\in\mathbb N$ depending only on $k$, such that for all pairs $(n,q)$, if
$q^{n^2}\ge N(k)$, then $P=M_{n,q}$. In particular, for a fixed $k$, the
conjecture implies that we have $P = M_{n,q}$ if either $q$ or $n$ is sufficiently large. The former part is solved by Theorem \ref{thm7} (first solved in \cite{KS23}) and the latter part is solved, in a sharper form than the solution provided by Theorem \ref{thm7} alone, by the next
4 theorems for every $q$ outside an explicitly computable finite set $E(k)$ of
powers of prime when either $p$ does not divide $k$ or only cyclic
matrices are considered (see Examples \ref{EX4} and \ref{EX5}).

\begin{theorem}\label{thm9}
Suppose $d\neq q-1$ $\big($thus $q-1\ge d(\ell+1)\big)$. If either $q-1\neq d(\ell+1)$ and $n>\ell(p-1)$ or $n>(\ell+1)(p-1)$, then the inclusion $\{A\in M_{n,q}|A {\rm\ is\ cyclic\ and\ }\Tr(A)\in\bbf_{p^{\ell}}\}\subset P$ holds. In particular, if $\ell=l$ (equivalently, $\mathbb F_q=d\mathbb P_{d,q}$, by Equation (\ref{EQ2})) and either $q-1\neq d(l+1)$ and $n>l(p-1)$ or $n>(l+1)(p-1)$, then $$\{A\in M_{n,q}|A\;\textup{is cyclic}\}\subset P.$$
\end{theorem}

\begin{theorem}\label{thm9b}
If the cyclic group $\bbf_p^{\times}\cap\mathbb P_{d,q}$ has at least 3 elements $\big($i.e., $\gcd\bigl((q-1)/d,p-1\bigr)\ge 3\big)$ and $n>p-1$, then we have an inclusion $$\{A\in M_{n,q}|A {\rm\ is\ cyclic\ and\ }\Tr(A)\in\bbf_p\}\subset P.$$
\end{theorem}

In the case of nilpotent cyclic matrices, the last two theorems complement \cite[Theorem 1.2]{KS23}.

The following theorems complement \cite[Corollary 8 and Theorem 5]{BM16} which imply that if $p\nmid k$ and $n>p$, then we have an inclusion $\{A\in M_{n,q}|A {\rm\ is\ cyclic\ and\ }\Tr(A)\in\bbf_p\}\subset P$. They are trivial if $k=1$ and hence one could either assume $k\ge 2$ or use the convention $\ln(0):=-\infty$.

\begin{theorem}\label{thm18}
If $p\ge 5$, $q>d^2$, $p\nmid k$, and the cyclic group $\bbf_p^{\times}\cap\mathbb P_{d,q}$ has at least 2 elements $\big($i.e., $\gcd\bigl((q-1)/d,p-1\bigr)\ge 2\big)$, then $$P=M_{n,q} \textrm{ for } n>\max\left\{\frac{p^2-3p+9}{2},\frac{4\ln{(k-1)}}{\ln{q}}\right\}.$$
\end{theorem}

\begin{theorem}\label{nlarge}
Suppose $q\ge 3d+1$ $($resp. $q=2d+1)$, $\ell=l$, and $p\nmid k$. Then $P=M_{n,q}$ if either $d>2$ and $n>\max\{\frac{d-1}{2},\frac{4\ln{(k-1)}}{\ln{q}}\}$ $($resp. and $n>\max\{2d-1,\frac{4\ln{(k-1)}}{\ln{q}}\})$ or $d=2$.

\end{theorem}

We also prove finer results when $q=p$ or $p\in\{2,3\}$; see Theorems \ref{thm10} and \ref{thm16}.

Section \ref{newS2} groups together different properties of the subsets $\mathbb P_{k,n,q}$. Section \ref{S2} presents properties of stabilizer subgroups and applications of double cosets of $S_{n,q}$ in $M_{n,q}$. Section \ref{S3} studies cyclic matrices and proves Theorem \ref{nonscalar}. All main results are proved in Section \ref{S4}. 

The cases $q\in\{2,3\}$ are special as we have an equality $2\mathbb P_{k,q}=\mathbb F_q$ for all $k$ if and only if $q\in\{2,3\}$. We study in detail the case $(n,q)=(2,2)$ or $(3,2)$ in Section \ref{S6} and the case $(n,q)=(2,3)$ in Section \ref{S7}. For example, the conclusions of Theorems \ref{thm2} (a) and \ref{thm5} (a) and (b) are not true if $q=2$ and $n=2$ (see Remark \ref{remark4}). Thus, from the point of view of statements that hold for all $n\ge 2$, Theorems \ref{thm2} and \ref{thm5} are optimal (in connection to Theorem \ref{thm2} (c) for $q=3$, see Example \ref{EX2}). We propose open problems in Section \ref{S8}. \\

\noindent{\bf Additional notations.} For $m,w\in\mathbb{N}$, let $0_{m,w}$ be the zero matrix of size $m\times w$; so $0_m=0_{m,m}$. For two square matrices $A$ and $B$ of sizes $m\times m$ and $w\times w$, let Diag$(A,B):=
\begin{pmatrix}
A & 0_{m\times w} \\
0_{w\times m} & B
\end{pmatrix}
 $. We define Diag$(A_1,\ldots,A_m)$ similarly as a block diagonal matrix with $A_1,\ldots,A_m$ on the diagonal, where $A_1,\ldots,A_m$ are square matrices (we always identify entries with $1\times 1$ matrices). All unlisted entries in a matrix are zero. Let $M^T$ denote the transpose of a matrix $M$. For $S\subset M_{n,q}$, let $S^T:=\{M^T\mid M\in S\}$. Let $\mathcal R^{\times}$ be the multiplicative group of units of $\mathcal R$; so $K^{\times}=K\setminus\{0\}$ and $M_{n,q}^{\times}=G_{n,q}$. Let $o(g)$ be the order of $g\in G_{n,q}$.  

\begin{center}
\sc{Acknowledgements}
\end{center}
We thank Michael Larsen for helpful conversations and correspondence.

\section{On powers}\label{newS2}

Each matrix in $M_n(K)$ is similar to a direct sum of generalized Jordan blocks $J_{f,m}$, with $m\in\mathbb N$ and $f(x)=x^r-\sum_{i=0}^{r-1} a_ix^i\in K[x]$ irreducible of degree $r:=\deg(f)\in\mathbb{N}$, where

\noindent
$J_{f,1}:=\begin{pmatrix}
    0&&&& a_0\\
    1&0&&&a_1\\
    &\ddots&\ddots&&\vdots\\
    &&&0&a_{r-2}\\
    &&&1&a_{r-1}
\end{pmatrix}$ and $J_{f,m}:=\begin{pmatrix}
    J_{f,1}&I_{r}&&& \\
    &J_{f,1}&I_r&&\\
    &&\ddots&\ddots&\\
    &&&J_{f,1}&I_{r}\\
    &&&&J_{f,1}
\end{pmatrix}$ for $m\ge 2$ has size $rm\times rm$. Note that $J_{f,m}$ is invertible if and only if $f(x)\neq x$, is semisimple if and only if $m=1$, and is split semisimple if and only if $r=1$.

Each generalized Jordan block $J_{f,m}$ is cyclic because, as a $K[J_{f,m}]$-module, $K^{mr}$ is isomorphic to $K[x]/(f^m)\cong K[J_{f,m}]$. More generally, a matrix $A\in M_n(K)$ is cyclic if and only if in a Jordan block decomposition of it, for each monic irreducible polynomial $f(x)\in K[x]$ there exists at most one generalized Jordan block $J_{f,m}$ for some $m\in\mathbb N$.

\begin{lemma}\label{order}
Let $G$ be a subgroup of $G_{n,q}$ and $g\in G$. If $\gcd(k,|G|)=1$, then $g=h^k$ for some $h\in G$. In particular, $g\in\mathbb{P}_{k,n,q}$.
\end{lemma}

\begin{proof}
As $o(g)$ divides $|G|$, $\gcd(k,o(g))=1$. Let $s,t\in\bbz$ be such that $sk+to(g)=1$. Thus $g=g^{sk+to(g)}=(g^s)^k$ and we can take $h:=g^s$.
\end{proof}

\begin{lemma}\label{embedding}
If $M_{m,q^r}=P_{k,m,q^r}$ and $\deg(f)=r$, then $J_{f,m}\in P_{k,mr,q}$.
\end{lemma}

\begin{proof}
There exists a noncanonical embedding of $M_{m,q^r}$ into $M_{mr,q}=\End_{\bbf_q}(\bbf_{q^{r}}^{m})$ by left multiplications, and we can regard $J_{f,m}\in M_{mr,q}$ as the element $J_{a,m}\in M_{m,q^r}$ for some $a\in\bbf_{q^r}$. Hence the lemma follows as $M_{m,q^r}=P_{k,m,q^r}$.
\end{proof}

\begin{prop}\label{exponent}
We have $\mathbb I_{n,q}=\mathbb P_{k,n,q}$ if and only if $e_{n,q}|k$. 
\end{prop}

\begin{proof}
Suppose $e_{n,q}|k$. Let $A\in M_{n,q}$. Considering a Jordan block decomposition, $A$ is similar to a matrix Diag($A_1,A_2$), where $A_1$ is invertible and $A_2$ is nilpotent. As $A_1^{e_{n,q}}$ is an identity matrix and $A_2^n$ is a zero matrix, it suffices to prove that $n\le e_{n,q}$. But $\bbf_{q^n}^{\times}$ is a subgroup of $G_{n,q}$ by the proof of Lemma \ref{embedding}, so $q^n-1$ divides $e_{n,q}$. Hence $n\le q^n-1\le e_{n,q}$.

If $\mathbb I_{n,q}=\mathbb P_{k,n,q}$, then for each $g\in G_{n,q}$, $g^k$ is an invertible idempotent and hence $g^k=I_n$, from which it follows that $e_{n,q}|k$.
\end{proof}

\begin{remark}\label{remark1}
We say that a matrix $A\in M_n(K)$ is \emph{of constant power ranks} if $A$ and $A^2$ have equal ranks (equivalently, for all $m\in\mathbb N$, $A$ and $A^m$ have equal ranks). Let $\mathbb P^{-}_{k,n,q}:=\{A^k|A\in M_{n,q}\ \textup{is of constant rank}\}$. For $k\in\mathbb N$ we have $\mathbb P^{-}_{k,n,q}\in\{\mathbb P^{-}_{s,n,q}|s\in\{1,\ldots,e_{n,q}\}\}$. The inclusion $\mathbb P^{-}_{k,n,q}\subset\mathbb P_{k,n,q}$ is an equality if and only if $k\ge n$. Hence for $k\ge n$ we have $P=2\mathbb P^{-}_{k,n,q}\in\{2\mathbb P^{-}_{s,n,q}|s\in\{1,\ldots,e_{n,q}\}\}$. 
Thus, for $k\in\mathbb N$ we have 
$$\mathbb P_{k,n,q}\in\{\mathbb P_{s,n,q}|s\in\{1,\ldots,n-1\}\}\cup\{\mathbb P^{-}_{s,n,q}|s\in\{1,\ldots,e_{n,q}\}\}.$$ 
Hence to study $P$ we can assume that $k\le e_{n,q}$.
\end{remark}

\begin{lemma}\label{kdivn}
We write $k=k_1k_2$ with $\gcd(k_2,q-1)=1$ and $k_1$ a divisor of a power of $q-1$ (equivalently, each prime divisor of $k_1$ divides $q-1$). If $k_1$ divides $n$, then $\{aI_n|a\in\mathbb F_q\}\subset\mathbb P_{k,n,q}\subset P_{k,n,q}$ and hence for $q>2$ we have $\{aI_n|a\in\mathbb F_q\}\subset Q_{k,n,q}$.
\end{lemma}
\begin{proof}
We can assume that $k_1=n$. For $a\in\mathbb F_q$, let $b\in\mathbb F_q$ be such that $a=b^{k_2}$. Let $A_b\in M_{n,q}$ be such that $\chi_{A_b}(x)=x^n-b$; so $bI_n=A_b^n$ by Cayley--Hamilton Theorem and hence $aI_n=A_b^k$.
\end{proof}

We have the following criterion on when a triangular matrix in $M_{n,q}$ is a $k$th power.

\begin{lemma}\label{lifting}
Let $A\in M_{n,q}$ be a triangular matrix. Suppose $p\nmid k$. If $A$ has at most one $0$ diagonal entry, and all other diagonal entries of it are in $\mathbb P_{d,q}\setminus\{0\}=\mathbb P_{k,q}\setminus\{0\}$, then $A\in\mathbb P_{k,n,q}$.
\end{lemma}
\begin{proof}
We prove this in the `upper' case, as the `lower' case follows by transposition. Let $\mathcal U$ be the subring of $M_{n,q}$ that consists of all upper triangular matrices. Let $\mathcal I$ be the two-sided ideal of $\mathcal U$ consisting of all upper triangular matrices with all diagonal entries being $0$. 

Assume $A\in\mathcal U$ is an upper triangular matrix that satisfies the requirement in the lemma. Then $A-x_1^k\in\mathcal I$ for some $x_1\in\mathcal U$ as the diagonal entries of $A$ are $k$th powers. We can assume that $x_1={\rm Diag}(\lambda_1,\ldots,\lambda_n)$ is such that for $i,j\in\{1,\ldots,n\}$ we have $\lambda_i^k=\lambda_j^k$  if and only if $\lambda_i=\lambda_j$. 

We prove by induction on $m\in\{1,\ldots,n\}$ that there exists $x_m\in\mathcal U$ such that $x_m-x_1\in\mathcal I$ and $A-x_m^k\in\mathcal I^m$. As $\mathcal I^n=0$, this will imply that $A=x_n^k$.

The statement is true for $m=1$ by our settings. Suppose it is true for some $m\in\{1,\ldots,n-1\}$, we want to show that it is true for $m+1$. By our assumption, there exists $x_m\in\mathcal U$ such that $A-x_m^k\in\mathcal I^m$. But its image $A-x_m^k+\mathcal I^{m+1}\in\mathcal I^m/\mathcal I^{m+1}$ may not be zero. Consider the map $\rho:\mathcal U\to\mathcal U$, $\rho(z):=(x_m+z)^k-x_m^k$; it leaves invariant the two-sided ideals of $\mathcal U$. Thus $\rho(\mathcal I^r)\subset\mathcal I^r$ for all $r\in\mathbb{N}$. For $z\in\mathcal I^m$, 
\[
\rho(z)\equiv\sum_{t=0}^{k-1}x_m^tzx_{m}^{k-1-t}\ {\rm (mod}\ \mathcal I^{2m}{\rm )}.
\] 
Hence $\rho$ induces a $\mathcal U/\mathcal I$-linear endomorphism $\bar{\rho}$ of $\mathcal I^m/\mathcal I^{m+1}$ by the rule
\[
\bar{\rho}(\bar{z}):=\sum_{t=0}^{k-1}\bar{x}^t\bar{z}\bar{x}^{k-1-t},
\]
where $\bar{x}:=x_1+\mathcal I\in\mathcal U/\mathcal I$ and we view $\mathcal I^m/\mathcal I^{m+1}$ as a $\mathcal U/\mathcal I$-bimodule. 

Notice that as a $\mathcal U/\mathcal I$-bimodule, $\mathcal I^{m}/\mathcal I^{m+1}$ decomposes into the direct sum of irreducible bimodules, which are just those 1-dimensional vector spaces $\bbf_q E_{ij}$ with $j=i+m$, $i,j\in\{1,\ldots,n\}$. The action of $\bar{\rho}$ on the vector space $\bbf_q E_{ij}$ is the multiplication by 
\[
\lambda_{i,j}:=\sum_{t=0}^{k-1}\lambda_i^t\lambda_j^{k-1-t}.
\]
If $\lambda_i=\lambda_j$, $\lambda_{i,j}=k\lambda_i^{k-1}\neq 0$ as $p\nmid k$ and $\lambda_i\neq 0$ as $i\neq j$ and there exists at most one $0$ diagonal entry. If $\lambda_i\neq\lambda_j$, $\lambda_{i,j}=\frac{\lambda_i^k-\lambda_j^k}{\lambda_i-\lambda_j}\neq 0$ by the choices made. As all $\lambda_{ij}$s are nonzero, $\bar{\rho}$ is an automorphism. Thus there exists $z\in\mathcal I^m$ such that $(x_m+z)^k-x_m^k+\mathcal I^{m+1}=A-x_m^k+\mathcal I^{m+1}$; so we can take $x_{m+1}:=x_m+z$. This ends the induction and the proof.
\end{proof}

\begin{cor}\label{generalLU}
Let $a_1,\ldots,a_n$ be the diagonal entries of $A\in M_{n,q}$. Suppose $p\nmid k$. If $a_i=b_i^k+c_i^k$ for $b_i,c_i\in\bbf_q$, $i\in\{1,\ldots,n\}$, and there is at most one $0$ in $b_1,\ldots,b_n$ and at most one $0$ in $c_1,\ldots,c_n$, then $A\in P$.
\end{cor}

\begin{proof}
 We write $A=L+U$ as the sum of a lower triangular matrix $L$ and an upper triangular matrix $U$, where the diagonal entries of $L$ (resp. $U$) are $b_1^k,\ldots,b_n^k$ (resp. $c_1^k,\ldots,c_n^k$). Then by Lemma \ref{lifting}, $L,U\in\mathbb{P}_{k,n,q}$, hence $A\in P$.
\end{proof}

For the proofs of Theorem \ref{thm7} and of the results from Theorems \ref{thm8} to \ref{nlarge}, we will need the following result proved in \cite[Chapter 6, Section 2, Corollary 1]{Jol73}.

\begin{theorem}\label{Joly}
Let $s\ge 2$ be an integer, let $a_1,\ldots,a_s,b\in\bbf_q^{\times}$, and let $k_1,\ldots,k_s\ge 1$ be integers. Denote by $N_s(b)$ the number of solutions of the equation $b=a_1x_1^{k_1}+\cdots+a_sx_s^{k_s}$ in $\bbf_q^{s}$. Then 
\[
|N_s(b)-q^{s-1}|\le \Delta q^{(s-1)/2} 
\]
where $\Delta:=(\delta_1-1)\cdots(\delta_s-1)$ with $\delta_i:=\gcd(k_i,q-1)$ for $i\in\{1,\ldots,s\}$,
and hence $N_s(b)>0$ if $q^{s-1}>\Delta^2$.
\end{theorem}

\begin{remark}\label{remark2}
Notice that $N_s(b)$ and $\Delta$ do not change if $k_i$ is replaced by $\delta_i$ for each $i\in\{1,\ldots,s\}$.
\end{remark}

\begin{corollary}\label{Jolyc}
Let $\eta_1,\eta_2\in\{0,\ldots,q-1\}$. For $i\in\{1,2\}$ let $E_i\subset\mathbb F_q$ be a subset with $\eta_i$ elements. If $q>(d-1)^4+2(\eta_1+\eta_2)d$, then for each $b\in\bbf_q^{\times}$, the equation $b=x^d+y^d$ has a solution in $(\mathbb F_q\setminus E_1)\times (\mathbb F_q\setminus E_2)$.
\end{corollary}

\begin{proof}
The number of solutions in $\mathbb F_q^2$ of the equation $b=x^d+y^d$ is $N_2(b)\ge\sqrt{q}[\sqrt{q}-(d-1)^2]$ by Theorem \ref{Joly}. For $q>(d-1)^2$, we have $\sqrt{q}/[\sqrt{q}+(d-1)^2]> 1/2$. Hence if $q-(d-1)^4> 2(\eta_1+\eta_2)d$, $N_2(b)>(\eta_1+\eta_2)d$, and the corollary follows from the fact that the number of solutions of the equation $b=x^d+y^d$ with either $x\in E_1$ or $y\in E_2$ is at most $(\eta_1+\eta_2)d$.
\end{proof}

For the next two corollaries, let the product decomposition $k=k_1k_2$ be as in Lemma \ref{kdivn}. 

\begin{corollary}\label{const}
If $n=n_1+\cdots+n_s$ with $s,n_1,\ldots,n_s\in\mathbb N$ such that 
 for $q\ge 4$ we have an inequality $n_i>\frac{4\ln{(d_{n_i}-1)}}{\ln{q}}$ for all $i\in\{1,\ldots,s\}$ with $d_{n_i}\ge 3$ and $k_1\nmid n_i$, then we have an inclusion $\{aI_n|a\in\mathbb F_q\}\subset P$. 
\end{corollary}

\begin{proof}
As $\bbf_q=2\mathbb P_{k,q}$ for $q\in\{2,3\}$, we can assume that $q\ge 4$. For each $i\in\{1,\ldots,s\}$, if $k_1|n_i$, then $aI_{n_i}\in P_{k,n_i,q}$ by Lemma \ref{kdivn}. If $k_1\nmid n_i$, by the proof of Lemma \ref{embedding}, we can regard $a I_{n_i}$ as the element $a\in\bbf_q\subset\bbf_{q^{n_i}}$, and we check that $a I_{n_i}\in P_{k,n_i,q}$. It suffices to show that $a\in 2\mathbb P_{k,q^{n_i}}$ if $a\neq 0$. As $a\neq 0$ and $n_i>\frac{4\ln{(d_{n_i}-1)}}{\ln{q}}$, we have $q^{n_i}>(d_{n_i}-1)^4$, hence $a\in 2\mathbb P_{k,q^{n_i}}$ by Corollary \ref{Jolyc} applied to $\eta_1=\eta_2=0$. 

We conclude that $aI_n=$Diag$(A_1,\ldots,A_s)\in P$.
\end{proof}

\begin{corollary}\label{constQ}
If $n=n_1+\cdots+n_s$ with $s,n_1,\ldots,n_s\in\mathbb N$ such that for all $i\in\{1,\ldots,s\}$ the number $p\frac{q^{n_i}-1}{d_{n_i}}$ is even and we have an inequality $n_i> \frac{\ln{\bigl((d_{n_i}-1)^4+4d_{n_i}\bigr)}}{\ln{q}}$ (resp. $n_i\ge 2$) if $(d_{n_i},q)\neq (1,2)$ and $k_1\nmid n_i$ (resp. if $(d_{n_i},q)=(1,2)$), then $\{aI_n|a\in\mathbb F_q\}\subset Q$. 
\end{corollary}

\begin{proof}
Let $a\in\mathbb F_q$. For each $i\in\{1,\ldots,s\}$, as in the proof of Corollary \ref{const}, it suffices to show that $a$ is a sum of two nonzero $k$th powers for $a\in\bbf_q\subset\bbf_{q^{n_i}}$. If $a\neq 0$ and $n_i>\frac{\ln{\bigl((d_{n_i}-1)^4+4d_{n_i}\bigr)}}{\ln{q}}$, then we have $q^{n_i}>(d_{n_i}-1)^4+4d_{n_i}$, hence $a$ is a sum of two nonzero $k$th powers by Corollary \ref{Jolyc}. If $a\neq 0$ and $(d_{n_i},q)=(1,2)$, then $n_i\ge 2$, hence $a$ is a sum of two nonzero $k$th powers as $2^{n_i}\ge 4$. Suppose $a=0$. If $p=2$, then $0=1+1$ is a sum of two nonzero $k$th powers. If $\frac{q^{n_i}-1}{d_{n_i}}$ is even, we have $d_{n_i}|\frac{q^{n_i}-1}{2}$, which implies that $-1$ is a $d_{n_i}$-th powers. Hence $0=1+(-1)$ is a sum of two nonzero $k$th powers. Notice that if $k_1| n_i$ and $q>2$, $0_{n_i}\in Q_{k,n_i,q}$ by Lemma \ref{kdivn}.

We conclude that $aI_n=$Diag$(A_1,\ldots,A_s)\in Q$.
\end{proof}

\section{Stabilizers and Double Cosets}\label{S2}

We consider a nonempty subset $S$ of $M_{n,q}$ stable under conjugation. 

\begin{lemma}\label{subgroup}
We have $L_S=R_S$. If $S=S^T$, then $L_S=L_S^{T}$.
\end{lemma}
\begin{proof}
For $g\in G_{n,q}$, we have $Sg=gSg^{-1}g=gS$. Hence $L_S=R_S$. If $S=S^T$, then for $g\in L_S$ we have $g^{T}S=g^TS^T=(Sg)^{T}=(gS)^{T}=S^{T}=S$, hence $L_S=L_S^{T}$. 
\end{proof}

\begin{lemma}\label{normal}
The group $L_S$ is a normal subgroup of $G_{n,q}$.
\end{lemma}
\begin{proof}
If $(g,h)\in G_{n,q}\times L_S$, then $ghg^{-1}S=ghSg^{-1}=gSg^{-1}=S$, hence $gh^{-1}g\in L_S$ and the lemma follows. 
\end{proof}

\begin{prop}\label{LP}
If $L_P\supseteq S_{n,q}$, then $P=M_{n,q}$.
\end{prop}
\begin{proof}
Define $B_{n,q}:=M_{n,q}\setminus G_{n,q}$. As $L_{P}=R_{P}\supseteq S_{n,q}$, to prove $B_{n,q}\subset P$, it suffices to prove that a representative of each class in 
 $S_{n,q}\backslash B_{n,q}/S_{n,q}$ belongs to $P$. But this is so as, by defining idempotents $D_0:=0_n$ and $D_{r}:={\rm Diag}(I_r,0_{n-r})$ for $r\in\{1,\ldots,n-1\}$, we have
\[
S_{n,q}\backslash B_{n,q}/S_{n,q}=\{[D_r]\mid r\in\{0,\ldots,n-1\}\}.
\]

Similarly, we have
\[
S_{n,q}\backslash G_{n,q}/S_{n,q}=\left\{[{\rm Diag}(I_{n-1},a)]\mid a\in\bbf_q^{\times}\right\}.
\]

For $b\in\bbf_q$, define $M_b:=$
$
\begin{pmatrix}
1 & b-1 \\
-1 & 1
\end{pmatrix}
=\left(
\begin{array}{cc}
1 & 0 \\
-1 & 0
\end{array}
\right)+\left(
\begin{array}{cc}
0 & b-1 \\
0 & 1
\end{array}
\right)$\\$=\left(
\begin{array}{cc}
1 & 0 \\
-1 & 0
\end{array}
\right)^k+\left(
\begin{array}{cc}
0 & b-1 \\
0 & 1
\end{array}
\right)^k\in P. 
$
As det$(M_b)=b$, we have the following class description $S_{n,q}\backslash G_{n,q}/S_{n,q}=\{[M_a]|a\in\bbf_q^{\times}$\} when $n=2$, and $S_{n,q}\backslash G_{n,q}/S_{n,q}=\{[{\rm Diag}(I_{n-2},M_a)]|a\in\bbf_q^{\times}$\} when $n>2$, in which case Diag$(I_{n-2},M_a)\in P$ for $a\in\bbf_q^{\times}$.

Therefore, $P\supset M_{n,q}=B_{n,q}\cup G_{n,q}$. Hence $P=M_{n,q}$.
\end{proof}

From the proof above, we notice that if we want to replace $P$ by $Q$ (resp. by $Q\cup\{0_n\}$), we only need to prove that the following two statements are true: 

\vskip 2mm
(i) For $r\in\{0,\ldots,n-1\}$ (resp. $r\in\{1,\ldots,n-1\}$) there exists a matrix in $Q$ of rank $r$. 

(ii) For $a\in \bbf_q^{\times}$ there exists a matrix in $Q$ of determinant $a$.

\begin{prop}\label{LQ}
We assume that $p$ does not divide $k$ and $L_Q\subset S_{n,q}$. Then $M_{n,q}=Q\cup\{0_n\}$. If moreover $-1\in\mathbb P_{d,q}$, then $Q=M_{n,q}$.
\end{prop}
\begin{proof}
If $-1\in\mathbb P_{d,q}=\mathbb P_{k,q}$, then $0=c^k+1^k$ for some $c\in\mathbb F_q^{\times}$, hence statement (i) for $Q$ follows from the facts that for $n=2$, the matrix $\begin{pmatrix}1&1\\1&1\end{pmatrix}=I_2^k+\begin{pmatrix}0&1\\1&0\end{pmatrix}^k$ if $k$ is odd and $\begin{pmatrix}2&2\\2&2\end{pmatrix}=\begin{pmatrix}1&0\\1& 1\end{pmatrix}+\begin{pmatrix}1 &1\\0&1\end{pmatrix}=\begin{pmatrix}1&0\\1/k&1\end{pmatrix}^k+\begin{pmatrix}1&1/k\\0&1\end{pmatrix}^k$ if $k$ is even (so $p$ is odd as it does not divide $k$) has rank $1$ and is in $Q$, and
$J_2=\begin{pmatrix}2& 1\\3 &2\end{pmatrix}=\begin{pmatrix}1 &1\\0& 1\end{pmatrix}+\begin{pmatrix}1 &0\\3&1\end{pmatrix}=\begin{pmatrix}1 &1/k\\0& 1\end{pmatrix}^k+\begin{pmatrix}1 &0\\3/k&1\end{pmatrix}^k$
has rank 2 and is in $Q$. 

If $-1\notin\mathbb P_{d,q}$, then $p$ is odd, and for $r\in\{1,\ldots,n-1\}$ let $J_r:=\begin{pmatrix}A_r\\B_r\end{pmatrix}$, where $A_r:=\begin{pmatrix}
    2I_{r-1} & 0_{(r-1)\times (n-r+1)}\end{pmatrix}$ is of size $(r-1)\times n$
and $B_r$ of size $(n-r+1)\times n$ has all entries $2$ (so if
$r=1$, then $J_r$ has all entries $2$). Clearly, $J_r$ has all diagonal
entries 2 and rank $r$ and we decompose $J_r=(I_n+L_r)+(I_n+U_r)$, where
$L_r$ is nilpotent lower triangular and $U_r$ is nilpotent upper triangular. The orders $o(I_n+L_r)$ and $o(I_n+U_r)$, being powers of $p$, are coprime to $k$, and thus $I_n+L_r,I_n+U_r\in\mathbb P_{k,n,q}$. So $J_r\in Q$ and statement (i) holds for $Q\cup\{0_n\}$. 

We check that statement (ii) is true. For each $b\in\bbf_q$, there exists a $2\times 2$ matrix $A_b\in Q$ of determinant $b$: $A_b:=\begin{pmatrix}2&1\\-b+4& 2\end{pmatrix}$ is equal to
\[
\begin{pmatrix}1&1\\0&1\end{pmatrix}+\begin{pmatrix}1&0\\-b+4&1\end{pmatrix}=\begin{pmatrix}1&1/k\\0&1\end{pmatrix}^k+\begin{pmatrix}1&0\\(-b+4)/k&1\end{pmatrix}^k.
\]
If $p\neq 2$, we can take Diag$(A_1,\ldots,A_1,A_a)\in Q$ if $n$ is even, and Diag$(2,A_1,\ldots,A_1,A_{a/2})$ if $n$ is odd, to have determinant $a\in\mathbb F_q$. If $p=2$, let $B_a$ be the invertible matrix $\begin{pmatrix}0&1&0\\0&0&1\\a&0&0\end{pmatrix}$. Let $L_a:=I_3+ak^{-1}E_{31}$ and
$U:=I_3+k^{-1}(E_{12}+E_{23})-k^{-3}\binom{k}{2}E_{13}$, where $E_{ij}$ denotes the $3\times 3$ matrix with entry $1$ in the $(i,j)$ position and with entry $0$ in all the other positions. Then $B_a=L_a^k+U^k$ and det$(B_a)=a$, and we can take Diag$(A_1,\ldots, A_1, A_a)\in Q$ if $n$ is even, and Diag$(A_1,\ldots,A_1,B_a)\in Q$ if $n$ is odd, to have determinant $a$. Hence statement (ii) holds.
\end{proof}

We recall the following classical result due to Jordan--Moore for $n=2$ and to Jordan--Dickson for $n>2$ (see \cite[Theorems 8.13 and 8.23]{Rot95}; see also \cite[Ch. VIII, Theorems 8.4 and 9.3]{Lan02} or \cite[Theorem 2.2.7]{GLS94}).

\begin{theorem}\label{simple}
The quotient group $S_{n,q}/Z_{n,q}$ is simple if and only if $(n,q)\notin\{(2,2),(2,3)\}$.
\end{theorem}

\begin{prop}\label{surjection}
Assume that $(n,q)$ is neither $(2,2)$ nor $(2,3)$. For the nonempty subset $S$ of $M_{n,q}$ stable under conjugation the following statements are equivalent:

(a) We have $L_{S}\cap (S_{n,q}\setminus Z_{n,q})\neq\varnothing$.

(b) The quotient map $\psi: L_S\cap S_{n,q}\rightarrow S_{n,q}/Z_{n,q}$ is a surjection.

(c) We have $L_S\supseteq S_{n,q}$.
\end{prop}

\begin{proof}
(c)$\Rightarrow$(a): this is clear.

(a)$\Rightarrow$(b): as $L_S$ is a normal subgroup of $G_{n,q}$ by Lemma \ref{normal}, $L_S\cap S_{n,q}$ is a normal subgroup of $S_{n,q}$. Hence the image $Im(\psi)$ of $\psi$ is normal. But by assumptions, $S_{n,q}/Z_{n,q}$ is simple (see Theorem \ref{simple}) and $Im(\psi)$ is not the identity. Hence $\psi$ is surjective.

(b)$\Rightarrow$(c): we consider the commutator map 
\[
\Gamma: S_{n,q}/Z_{n,q}\times S_{n,q}/Z_{n,q}\to S_{n,q},\ (\bar{g},\bar{h})\mapsto[g,h]={g}{h}{g}^{-1}{h}^{-1}.
\]
As $\psi$ is surjective, we deduce that
\[
L_S\cap S_{n,q}\supseteq[L_S\cap S_{n,q},L_S\cap S_{n,q}]=[S_{n,q},S_{n,q}]=S_{n,q}
\]
(e.g., see \cite[Ch. VIII, Theorems 8.3 and 9.2]{Lan02} for last equality).
\end{proof}

\begin{example}\label{EX1}
The nonempty subset $\Pi_{n,q}\subset M_{n,q}$ is stable under conjugation, transposition, and the set involution of $M_{q,n}$ defined by the rule $A\mapsto 2I_n-A$. Recall that $\Pi_{n,q}\subsetneq M_{n,q}$ by \cite[Lemma 3]{HT88}.

We first assume that $q\ge 4$. Let $a\in\mathbb F_q^{\times}\setminus\{1,2\}$. We show that the assumption that $aI_n\in\Pi_{n,q}$ leads to a contradiction. This assumption implies
that there exist $g,h\in G_{q,n}$ and $s,t\in\{0,\ldots,n\}$ such that
$aI_n=gD_sg^{-1}+hD_th^{-1}$ with $D_r={\rm Diag}(1,\ldots,1,0,\ldots,0)\in\mathbb I_{n,q}$ of rank $r\in\{0,\ldots,n\}$. Conjugating with $g^{-1}$ we can assume that $g=I_n$, hence $aI_n-D_s=hD_th^{-1}\in\mathbb I_{n,q}$, which is not possible as the
eigenvalues of $aI_n-D_s$ belong to the set $\{a,a-1\}$ and we have $\{a,a-1\}\cap\{0,1\}=\varnothing$. 

If $p=3$, then a similar argument gives that $2I_n+E\notin\Pi_{n,q}$, where $E$ is a nonzero upper triangular matrix with $E^2=0_n$; there exist such matrices $E$ with the property that $I_n+2^{-1}E\in\Pi_{n,q}$.

The proof of Proposition \ref{LP} gives: if $L_{\Pi_{n,q}}\supseteq S_{n,q}$, then $\Pi_{n,q}=M_{n,q}$. Hence we have $L_{\Pi_{n,q}}\cap S_{n,q}\neq S_{n,q}$. From this and Proposition \ref{surjection} it follows that for $(n,q)\notin\{(2,2),(2,3)\}$ we have $L_{\Pi_{n,q}}\cap S_{n,q}\subset Z_{n,q}$, and based on the last two paragraphs we conclude that $L_{\Pi_{n,q}}\cap S_{n,q}=\{I_n\}$.
\end{example}

\begin{proof}[Proof of Theorem \ref{thm1}]
If $g\in L_P$, from Propositions \ref{surjection} and \ref{LP} it follows that $P=M_{n,q}$. On the other hand, if $P=M_{n,q}$, then $g\in G_{n,q}=L_P$. Hence part (a) holds. The proofs of parts (b) and (c) are the same but with Proposition \ref{LQ} replaced by Proposition \ref{LP}.
\end{proof}

\section{Diagonal entries in conjugacy classes}\label{S3}

\begin{lemma}\label{cyclic}
For $A\in M_n(\mathcal R)$, with $\mathcal R$ a commutative ring, the following statements are equivalent:

(1) The matrix $A$ is a cyclic matrix.

(2) If $(u_1,\ldots,u_n)\in\mathcal R^n$ is such that $\sum_{i=1}^{n}u_{i}=\Tr(A)$, then there exist $a_0,\ldots,a_{n-2}\in\mathcal R$ such that $A$ is similar to the matrix
\begin{equation}\label{EQ3}
\mathcal A:=\begin{pmatrix}
 u_1 &&&& a_0\\
 1 & u_2 &&& a_1\\
& \ddots &\ddots && \vdots\\
&& 1 & u_{n-1} & a_{n-2}\\
&&& 1 & u_n
\end{pmatrix}.
\end{equation}
\end{lemma}

\begin{proof}
$(1)\Rightarrow(2)$: Let $v\in\mathcal R^n$ be such that $\{v,Av,\ldots, A^{n-1}v\}$ is an $\mathcal R$-basis of $\mathcal R^n$.
We will construct recursively a new $\mathcal R$-basis $\{w_1,\ldots,w_n\}$ of $\mathcal R^n$ that has the following two properties:

(a) for each $i\in\{1,\ldots,n\}$, we have $\sum_{j=1}^i \mathcal Rw_j=\sum_{j=0}^{i-1} \mathcal RA^jv$ and
$w_i-A^{i-1}v\in\sum_{j=0}^{i-2} \mathcal RA^jv$ (so for $i=1$ we have $w_1=v$);

(b) the matrix representation $\mathcal A$ with respect to this $\mathcal R$-basis of the $\mathcal R$-linear map
$ A:\mathcal R^n\to\mathcal R^n$ is as in Equation (\ref{EQ3}).

Define $w_1:=v$. Once $w_1,\ldots,w_i$ are constructed and $i<n$, choose $w_{i+1}$ such that $A(w_i)=w_{i+1}+u_iw_i$. So $w_{i+1}:=A(w_i)-u_iw_i$. It is clear that all desired properties hold.

$(2)\Rightarrow(1)$: Take $u_1=\cdots=u_{n-1}=0$. Then $A$ is similar to a matrix $\mathcal A$ as in Equation (\ref{EQ3}). Let $v_1:=(1,0,\ldots,0)^T$. It is easy to see that $\{v_1,\mathcal Av_1,\ldots,\mathcal A^{n-1}v_1\}$ is an $\mathcal R$-basis of $\mathcal R^n$, hence $A$ is cyclic.
\end{proof}

\begin{cor}\label{cyclicLU}
We assume $A\in M_{n}(\mathcal{R})$ is cyclic. Then for every vector $(e_1,u_1,\ldots,e_n,u_n)\in\mathcal{R}^{2n}$ such that $\Tr(A)=\sum_{i=1}^n(e_i+u_i)$, there exist $a_0,\ldots,a_{n-2}\in\mathcal{R}$ such that $A$ is similar to the sum $\mathcal{A}=B+C$:
\begin{equation}
 \mathcal{A}:=\begin{pmatrix}
    e_1 &&&&&\\
    0 & e_2 &&&&\\
    &1& e_3 &&&\\
    &&0&e_4&&\\
    &&&\ddots&\ddots&\\
    &&&&1&e_n
\end{pmatrix}+\begin{pmatrix}
    u_1 &&&&&a_0\\
    1 & u_2 &&&&a_1\\
    &0&u_3&&&a_2\\
    &&1&u_4&&a_3\\
    &&&\ddots&\ddots&\vdots\\
    &&&&0&u_n
\end{pmatrix}
\end{equation}
if $n$ is odd;
\begin{equation}
 \mathcal{A}:=\begin{pmatrix}
    e_1 &&&&&\\
    1 & e_2 &&&&\\
    &0& e_3 &&&\\
    &&1&e_4&&\\
    &&&\ddots&\ddots&\\
    &&&&1&e_n
\end{pmatrix}+\begin{pmatrix}
    u_1 &&&&&a_0\\
    0 & u_2 &&&&a_1\\
    &1&u_3&&&a_2\\
    &&0&u_4&&a_3\\
    &&&\ddots&\ddots&\vdots\\
    &&&&0&u_n
\end{pmatrix}
\end{equation}
if $n$ is even.
\end{cor}

\begin{example}\label{decompex}
If $n$ is odd (resp. even), then $B$ is split semisimple provided $e_{2i}\neq e_{2i+1}$ (resp. $e_{2i-1}\neq e_{2i}$) for all $i\in\{1,\ldots,\lceil (n-1)/2\rceil\}$. If $n$ is odd (resp. even), then $C$ is split semisimple provided $u_{2i-1}\neq u_{2i}$ (resp. $u_{2i}\neq u_{2i+1}$) for all $i\in\{1,\ldots,\lfloor (n-1)/2\rfloor\}$ and $u_n\neq u_j$ for all $j\in\{1,\ldots,n-1\}$. To see this, one can look at the rank of $B-e_iI_n$ or $C-u_iI_n$ for $i\in\{1,\ldots,n\}$ and deduce that the sum of the dimensions of the eigenspaces of either $B$ or $C$ equals to $n$. E.g., for $p\ge 5$ and 
\[
\mathcal{A}=\begin{pmatrix}
    1 &&&&&\\
    1 & 0 &&&&\\
    &0& 2 &&&\\
    &&1& 0&&\\
    &&&0&2&\\
    &&&&1&0
\end{pmatrix}+\begin{pmatrix}
    0 &&&&&a_0\\
    0 & 1 &&&&a_1\\
    &1&0&&&a_2\\
    &&0&2&&a_3\\
    &&&1&0&a_4\\
    &&&&0&3
\end{pmatrix},
\]the dimensions of the eigenspaces for the eigenvalues $0,1,2$ of $B$ are $3,1,2$ (respectively) and the dimensions of the eigenspaces for the eigenvalues $0,1,2,3$ of $C$ are $3,1,1,1$ (respectively), and they add up to 6.
\end{example} 

\begin{cor}\label{diagonal}
Let $K$ be a field with at least 3 elements. Then each $A\in M_n(K)$ is similar to Diag$(A_1,\ldots,A_s)$, where $s\in \mathbb{N}$ and each $A_i$ is a square matrix which is either $0$ or has all diagonal entries nonzero. 
\end{cor}

\begin{proof}
As each generalized Jordan block is cyclic, this follows from Lemma \ref{cyclic}. Notice that, as $K$ has at least 3 elements, for each element $a\in K$, there exist $u_1,\ldots,u_n\in K^{\times}$ such that $u_1+\cdots+u_n=a$.
\end{proof}

\begin{defn}\label{defqcyc}   
Let $\mathcal R$ be a commutative ring. We say that a matrix $A\in M_n(\mathcal R)$ is \emph{quasi-cyclic} if for every vector $(a_1,\ldots,a_n)\in\mathcal R^n$ such that $\sum_{i=1}^n a_i=\Tr(A)$, $A$ is similar to a matrix whose diagonal entries are $a_1,\ldots,a_n$.
\end{defn}

\begin{lemma}\label{q1}
Let $A\in M_n(\mathcal R)$ and $b\in\mathcal R$. Then $A$ is quasi-cyclic if and only if $A-bI_n$ is quasi-cyclic. 
\end{lemma}
\begin{proof}
As $A=(A-bI_n)+bI_n$, it suffices to prove the `only if' part. So, assuming $A$ is quasi-cyclic, to check that $A-bI_n$ is cyclic, let $(u_1,\ldots,u_n)\in\mathcal R^n$ with $u_1+\cdots+u_n=\Tr(A-bI_n)$. Let $g\in GL_n(\mathcal R)$ be such that $gAg^{-1}$ has diagonal entries $u_1+b,\ldots,u_n+b$. As $g(A-bI_n)g^{-1}$ has diagonal entries $u_1,\ldots,u_n$, $A-bI_n$ is quasi-cyclic.
\end{proof}

\begin{lemma}\label{q2}
Let $a,b\in\mathcal R$ be such that $a-b\in\mathcal R^{\times}$. Then the two matrices $D_1:={\rm Diag}(a,bI_m)$ and $D_2:={\rm Diag}\left(\begin{pmatrix}
    b & 1\\
    0 & b
\end{pmatrix},bI_m\right)$ are quasi-cyclic for all $m\in\mathbb N$.
\end{lemma}

\begin{proof}
By Lemma \ref{q1} we can assume that $b=0$; so $a\in\mathcal R^{\times}$. Let $i\in\{1,2\}$. For $w:=m+i$, let $(u_1,\ldots,u_w)\in\mathcal R^w$ be such that $\sum_{i=1}^w u_i$ is $a$ if $i=1$ and is $0$ if $i=2$. The product matrix 
$$M:=(1 u_2 \cdots u_w)^T\cdot (u_1 1 \cdots 1)\in M_w(\mathcal R)$$ 
has diagonal entries $u_1,\ldots,u_w$ and its matrix representation with respect to the $\mathcal R$-basis $\{M(v_w), v_w-a^{-1}M(v_w),v_3-v_2,\ldots,v_w-v_{w-1}\}$ of $\mathcal R^w$ is $D_i$, where $\{v_1,\ldots,v_w\}$ is the standard $\mathcal R$-basis of $\mathcal R^w$. Hence $M$ and $D_i$ are similar. Thus $D_i$ is quasi-cyclic. 
\end{proof}

\begin{lemma}\label{q3}
If $n\in\{2,3\}$, then each nonscalar matrix $A\in M_n(K)$ is quasi-cyclic. 
\end{lemma}
\begin{proof}
By Lemma \ref{cyclic} we can assume that $A$ is not cyclic. Considering a Jordan block decomposition of $A$ it follows that $n=3$ and that $A$ is similar to either Diag$(a,b,b)$ or Diag $\left(\begin{pmatrix}
    a & 1 \\
    0 & a
\end{pmatrix},a\right)$, where $a,b\in K$ are distinct; from Lemma \ref{q2} it follows that $A$ is quasi-cyclic.
\end{proof}

Let $A\in M_n(K)$. For a nonempty subset $\tau\subset\{1,2,\ldots, n\}$ with $m$ elements, let $A|_{\tau}\in M_{m}(K)$ be the submatrix of $A$ formed from the $i$-th rows and the $j$-th columns of $A$ with $i,j\in\tau$.

The following result was first proved in \cite[Theorem 2 and comments after it]{Fil69}, and is reproved independently here in a way that partially applies to arbitrary commutative rings $\mathcal R$.

\begin{theorem}\label{nonscalar}
Suppose $A\in M_n(K)$. Then $A$ is a nonscalar matrix if and only if $A$ is quasi-cyclic.
\end{theorem}

\begin{proof}
We prove by induction on $n\ge 2$ that each nonscalar $A\in M_n(K)$ is quasi-cyclic. The base of the induction for $n\le 3$ holds by Lemma \ref{q3}. For $n\ge 4$, the passage from $\le n-1$ to $n$ goes as follows. Up to similarity, based on Lemma \ref{cyclic} we can write $A=$Diag$(Q_1,...,Q_s,S_1,...,S_t)$ with each $Q_i$ quasi-cyclic and every $S_j=a_jI_{m_j}$ scalar. We will use a second induction on
$w:=s+t$. If $w=1$, then, as $A$ is nonscalar, $A=Q_1$ is quasi-cyclic. Assume
that for some $m\ge 2$, each $A$ is quasi-cyclic if $w<m$, and we consider
the case $w=m$. If $t>0$, we can assume that $a_1,\ldots,a_t$ are
pairwise distinct. Thus, by the induction on $n$, we can assume that $s,t\le 2$, that
for $s>0$ we have $t\le 1$, and that for $t>0$ we have $s\le 1$. The case
$(s,t)=(0,2)$ gets reduced to the case $s=t=1$ by the induction on $n$.
Thus, as $w=m\ge 2$, we can assume that $m=2$ and $(s,t)\in\{(2,0),(1,1)\}$. 

Let $(u_1,\ldots,u_n)\in K^n$ be such that $\sum_{i=1}^n u_i=\Tr(A)$.

If $A=$Diag$(Q_1,Q_2)$, for $i\in\{1,2\}$ let $n_i\times n_i$ be the size of $Q_i$; so $n_i\ge 2$ and $n=n_1+n_2$. As $Q_1$ is quasi-cyclic, it is similar to a matrix $P_1$ with diagonal entries $u_1,\ldots,u_{n_1-1},\Tr(Q_1)-\sum_{i=1}^{n_1-1} u_i$. Hence $A$ is similar to $B:=$Diag$(P_1,Q_2)$. As $Q_2$ is a nonscalar matrix, so is $B|_{\{n_1,n_1+1,\ldots,n\}}$. So $B|_{\{n_1,n_1+1,\ldots,n\}}$ is similar to a matrix with diagonal entries $u_{n_1},\ldots,u_n$. Hence $B$ is similar to a matrix with diagonal entries $u_1,\ldots,u_n$. Thus $A$ is quasi-cyclic.

If $A=$Diag$(Q_1,S_1)$, by the induction on $n$ and Lemma \ref{q1}, we can assume
that $S_1=0_1$. As $Q_1$ is quasi-cyclic, it is similar to a matrix $P_1$ with last and first diagonal entries equal to $1$ and $u_0\in K\setminus\{1-u_n\}$. So $A$ is similar to $B:=$Diag$(P_1,0)$. As $1\neq 0$, $B|_{\{n-1,n\}}$ is nonscalar, hence similar to a matrix with diagonal entries $1-u_n,u_n$. So $B$ is similar to a matrix $C$ with the last $2$ entries equal to $1-u_n,u_n$ and  same first $n-2$ entries. As $u_0$ and $1-u_n$ are distinct diagonal entries of $C|_{\{1,\ldots,n-1\}}$, this matrix is nonscalar, hence by the induction hypothesis on $n$, it is quasi-cyclic and thus similar to a matrix with diagonal entries $u_1,\ldots,u_{n-1}$. Hence $C$ is similar to a matrix with diagonal entries $u_1,\ldots,u_n$. Thus $A$ is quasi-cyclic. This ends both inductions.

As scalar matrices are not quasi-cyclic, the theorem holds.\end{proof}

\begin{example}\label{EX3} 
Let $A\in M_{n,q}$ be nonscalar. If $\Tr(A)-2n\in\{-2,-1,0\}$ (e.g., this holds if $q\in\{2,3\}$) and $p\nmid k$, then from Corollary \ref{generalLU} applied with $B$ and $C$ having all the first $n-1$ diagonal entries $1$ and the last diagonal entries in $\{0,1\}$ and from Theorem \ref{nonscalar} it follows that $A\in P$. 
\end{example}

As for each monic polynomial $f(x)\in K[x]$ of degree $n$ there exists a nonscalar matrix $A\in M_n(K)$ such that $\mu_A(x)=f(x)$, Theorem \ref{nonscalar} reobtains and refines \cite[Theorem 2.1]{FL58}. It is not clear how other proofs and refinements of \cite[Theorem 2.1]{FL58}, such as \cite[Theorems 1 to 3]{DeO73} and \cite{Her83}, could help in the study $P$. 

\section{Proofs of the main theorems}\label{S4}

\begin{proof}[Proof of Theorem \ref{thm2} (a)]
We include a proof that works for $q\ge 3$ and $k$ not divisible by $p$; so, as $d=1$, we have $\gcd\bigl(k,p(q-1)\bigr)=1$. Let $A\in M_{n,q}$. We write $A=L+U$ with $L,U\in G_{n,q}$ such that $L$ is lower triangular and $U$ is upper triangular as follows. As $q>2$, there exists $b\in\bbf_q^{\times}\backslash\{1\}$. For $i\in\{1,\ldots,n\}$, let $a_i$ be the $(i,i)$ entry of $A$. If $a_i\neq 1$, let $u_i:=1$. If $a_i=1$, let $u_i:=b$. Then we take $L,U\in G_{n,q}$ such that $a_1-u_1,\ldots,a_n-u_n$ are the diagonal entries of $L$ and $u_1,\ldots,u_n$ are the diagonal entries of $U$.

As $L$ belongs to the subgroup of order $(q-1)^{n}q^{n(n-1)/2}$ of $G_{n,q}$ that consists of invertible lower triangular matrices, we have $L\in\mathbb P_{k,n,q}$ by Lemma \ref{order} as $\gcd(k,p(q-1))=1$. Similarly, we argue that $U\in\mathbb P_{k,n,q}$. Therefore $Q=M_{n,q}$.
\end{proof}

\begin{theorem}\label{thm3}
If $q=2$ and $k$ is odd, then $P=M_{n,q}$.
\end{theorem}

\begin{proof}
It suffices to show that every generalized Jordan block $J_{f,m}$ is a sum of two $k$th powers. The diagonal entries of $J_{f,1}$ are $0,0,\ldots,0,a_{r-1}$, where $r=\deg(f)$ and $a_{r-1}\in\{0,1\}$. We decompose
\[
J_{f,1}=\begin{pmatrix}
     1&&&0 \\
    1&1&&0\\
    &\ddots&\ddots&\vdots\\
    &&1&(a_{r-1}-1)  
\end{pmatrix}+\begin{pmatrix}
    1&&&& a_0\\
    &1&&&a_1\\
    &&\ddots&&\vdots\\
    &&&&1
\end{pmatrix}
\]
and we denote the first and the second terms of the sum by $L$ and $U$ (respectively). Then $o(U)$ is a a power of $2$. As $\gcd(k,2)=1$, $U\in\mathbb P_{k,r,q}$ (see Lemma \ref{order}). If $a_{r-1}=0$, then $L\in\mathbb P_{k,r,q}$ by the same argument. If $a_{r-1}=1$, then $L$ is similar to Diag$(J_{x-1,r-1},0)\in\mathbb P_{k,r,q}$. Hence $J_{f,1}\in P$. For $m>1$, $J_{f,m}=L_m+U_m$, where $L_m:=$Diag$(L,\ldots,L)$ and $U_m:=\begin{pmatrix}
    U&I_r&&\\
    &&\ddots&\ddots&\\
    &&&U&I_{r}\\
    &&&&U
\end{pmatrix}$, and we similarly argue that $J_{f,m}\in P$. 
\end{proof}

\begin{theorem}\label{thm4}
If $q=2$ and $\gcd(k,6)=1$, then $Q=M_{n,q}$.    
\end{theorem}

\begin{proof}
The proof is the same as above, except that, based on \cite[Theorem 12 or Corollary 13]{Hen74} or on \cite[Theorem]{Lor87}, we write the block matrix $\begin{pmatrix}
    0& a_{r-2}\\
    1& a_{r-1}
\end{pmatrix}$ in $J_{f,1}$ as a sum of two invertible matrices. Now, as $\gcd(k,6)=1$ and $|G_{2,2}|=6$, those two invertible matrices are $k$th powers. Hence $Q=M_{n,q}$.
\end{proof}

\begin{remark}\label{remark3}
Each matrix in $M_{n,2}$ is the sum of an idempotent and of a unipotent matrix (see \cite[Theorem 3]{BCDM13}) and based on this one easily gets a new proof of Theorem \ref{thm3}.

Similarly, if $d=1$ and $p$ is odd and does not divide $k$, then the identity $P=M_{n,q}$, which is a weaker form of Theorem \ref{thm2} (b), also follows easily from \cite[Theorem 1]{Bre18}.
\end{remark}

\begin{theorem}\label{thm5}
If $q\ge3$, then each matrix in $M_{n,q}$ is a sum of

(a) two invertible cyclic matrices; 

(b) two invertible semisimple matrices.
\end{theorem}

\begin{proof} (a) As in the proof of \cite[Theorem]{Lor87}, it suffices to check that the number $c_{n,q}$ of invertible cyclic matrices in $G_{n,q}$ is more than $|M_{n,q}|/2$. But $c_{n,q}\ge|G_{n,q}|\big(1-1/q(q^2-1)\big)$ by \cite[Theorem 3.1 and Remark 3.2]{NP00}. Hence it suffices to show
\[
2\left(1-\frac{1}{q(q^2-1)}\right)(q^n-1)(q^n-q)\cdots(q^n-q^{n-1})>q^{n^2},
\]
which is equivalent to 
\begin{equation}\label{EQ4}
2\left(1-\frac{1}{q(q^2-1)}\right)(1-q^{-n})(1-q^{-n+1})\cdots(1-q^{-1})>1.
\end{equation}
But $(1-q^{-n})(1-q^{-n+1})\cdots(1-q^{-1})>1-q^{-1}-q^{-2}$ by the pentagonal number theorem, hence the Inequality (\ref{EQ4}) holds as for the increasing function $F:[3,\infty)\to\mathbb R$, $F(x)=2\left(1-\frac{1}{x(x^2-1)}\right)(1-x^{-1}-x^{-2})$, we have $F(3)=\frac{23}{12}(1-\frac{1}{3}-\frac{1}{9})=\frac{115}{108}>1$.

(b) It suffices to prove that each Jordan block $J_{f,m}\in M_{mr,q}$ is a sum of two invertible semisimple matrices. We can further assume that $f(x)=x-a$ with $a\in\bbf_q$ by the embedding part of the proof of Lemma \ref{embedding} as it maps (invertible) semisimple elements to (invertible) semisimple elements. We will decompose
\[
J_{x-a,m}=\begin{pmatrix}
    a&1&&& \\
    &a&1&&\\
    &&\ddots&\ddots&\\
    &&&a&1\\
    &&&&a\\
\end{pmatrix}=B+C,
\]
with $B,C\in G_{m,q}$ semisimple.

If either $q\ge 4$ or $q=3$ with $a=0$, there exist distinct $b_1,b_2\in\bbf_q^{\times}$ such that $a=b_1+b_2$. Thus we can take 
\[
B:=\begin{pmatrix}
    b_1&1&&&& \\
    &b_2&0&&&\\
    &&b_1&1&&\\
    &&&\ddots&\ddots&\\
    &&&&b_i&2-i\\
    &&&&&b_{3-i}\\
 
\end{pmatrix}, \; C:=\begin{pmatrix}
    b_2&0&&&& \\
    &b_1&1&&&\\
    &&b_2&0&&\\
    &&&\ddots&\ddots&\\
    &&&&b_{3-i}&i-1\\
    &&&&&b_i\\
\end{pmatrix},
\]
where $i\in\{1,2\}$ is such that $m-i$ is odd, $b_1$ and $b_2$ alternate on the main diagonals and $0$ and $1$ alternate on the diagonals of entries $(j-1,j)$ with $j\in\{2,\ldots,m\}$.

For $q=3$ and $a\in\{1,2\}$, the case of $J_{x-2,m}$ will follow from the case of $J_{x-1,m}$ by adding $I_m$ to $B$. For $a=1$ let $B_2:=\begin{pmatrix}
0 & 1\\ 1 & 2
\end{pmatrix}$, $C_2:=\begin{pmatrix}
1 & 0\\ 2 & 2
\end{pmatrix}$, and $F_2:=\begin{pmatrix}
0 & 0\\ 1 & 0
\end{pmatrix}$, and we notice that $B_2,C_2$ are semisimple as $\chi_{B_2}(x)=x^2+x+2\in\mathbb F_3[x]$ is irreducible and $\chi_{C_2}=(x-1)(x-2)$. 

If $m$ is even, we take $B:=$Diag$(B_2,\ldots,B_2)$ and 
\[
C:=\begin{pmatrix}
C_2 & F_2 & &\\
& \ddots & \ddots&\\
& & C_2 & F_2\\
& & & C_2
\end{pmatrix}.
\]
The even columns of $C-I_m$ or $C_2-I_m$ are 0 so the dimensions of the eigenspaces of $C$ of eigenvalues $1$ and $2$ add up to at least $m/2+m/2=m$ and thus to exactly $m$; thus $C$ is invertible split semisimple.

If $m$ is odd, let $B:=$Diag$(B_2,\ldots,B_2,2)$ and 
\[
C:=\begin{pmatrix}
C_2 & F_2 & &\\
& \ddots & \ddots&\\
& & C_2 & F_3\\
& & & C_3
\end{pmatrix},
\]
where $C_3:=\begin{pmatrix}
1 & 0 & 0\\ 2 & 2 & 1 \\ 0 & 0 & 2
\end{pmatrix}$ and $F_3:=(F_2\ 0_{2\times 1})$, and as above we argue that $C$ is invertible semisimple. 

We conclude that for all $q\ge 3$ and $a\in\mathbb F_q$, $J_{x-a,m}$ is a sum of two invertible semisimple matrices in $M_{m,q}$.
\end{proof}

\begin{proof}[Proof of Theorem \ref{thm2} (b) and (c)]
As $d=1$, $\bbf_q=\mathbb P_{k,q}$ and thus each split semisimple matrix over $\mathbb F_q$ is a $k$th power.

As each generalized Jordan block $J_{f,m}$ is cyclic, hence similar to a matrix $\mathcal A$ as in Equation (\ref{EQ3}) but of size $s\times s$ with $s:=m\deg(f)$, for $q\ge 3$ (resp. $q\ge 4$) it suffices to write $\mathcal A$ as a sum of two (resp. two invertible) split semisimple $s\times s$ matrices. Let $t:=\Tr(\mathcal{A})$.

We can assume $s\ge 2$. With $a,b\in\mathbb F_q$ distinct, we consider two cases.

Case 1: $s$ is odd. We write $\mathcal A=B+C$ by Corollary \ref{cyclicLU} using $(a,-a,b,-b,\ldots,a,-a,b,-b,c,t-c)\in\bbf_q^{2s}$ with $c\in\bbf_q\setminus\mathcal P_1$, where $\mathcal P_1:=\{b,t+a,t+b\}$. We have $\mu_B(x)=(x-a)(x-b)(x-c)$ if $c\neq a$ and $\mu_B(x)=(x-a)(x-b)$ if $c=a$, and $\mu_C(x)=(x+a)(x+b)(x-t+c)$, and $B$ and $C$ are split semisimple matrices over $\bbf_q$ by Example \ref{decompex} if such a $c$ exists for a good choice of the pair $(a,b)$. Clearly, $c$ exists if $q\ge 3$ provided for $q=3$ and $t\in\{1,2\}$ we take $b=a+t$.

Case 2: $s$ is even. We write $\mathcal A=B+C$ by Corollary \ref{cyclicLU} using $(a,-a,b,-b,\ldots,a,-a,b,-b,a,-a,c,t-c)\in\bbf_q^{2s}$ with $c\in\mathbb F_q\setminus\mathcal P_0$, where $\mathcal P_0:=\{a,t+a,t+b\}$. This case is similar to the Case 1 with the roles of $a$ and $b$ interchanged.

From the above two cases it follows that $P^{ss}=M_{n,q}$ if $q\ge 3$.

Assume now that $q\ge 4$. We take $a,b\in\mathbb F_q^{\times}$. Then $B$ and $C$ are also invertible provided $c\notin \mathcal Q_{\epsilon}:=\mathcal P_{\epsilon}\cup\{0,t\}$, where $\epsilon\in\{0,1\}$ is such that $s-\epsilon$ is even. Clearly, such a $c$ exists if $q\ge 7$ regardless of what the pair $(a,b)$ is, hence we can assume that $q\in\{4,5\}$. If $\epsilon=1$, then $\mathcal Q_{\epsilon}$ has three elements provided for $t\neq 0$ we take $b=t$ if $p=2$ and $(a,b)=(-2t,-t)$ if $p=5$. The case $\epsilon=0$ gets reduced to the case $\epsilon=1$ by interchanging the roles of $a$ and $b$. Hence $c$ exists for a good choice of the pair $(a,b)$ even if $q\in\{4,5\}$. 

From the last paragraph it follows that $Q^{ss}=M_{n,q}$ if $q\ge 4$.
\end{proof}

\begin{example}\label{EX2}
We assume that $(n,q)=(2,3)$. Each invertible split semisimple matrix in $M_{2,3}$ either has trace 0 or is scalar. Let $A:=\begin{pmatrix}
    2 & 1 \\ 0 & 2 
\end{pmatrix}$. As $\Tr(A)=1$, if we can decompose $A=B+C$ with $B$ and $C$ invertible split semisimple, then the intersection $\{B,C\}\cap\{I_2,2I_2\}$ is nonempty, and we reach a contradiction as $A-I_2$ and $A-2I_2$ are not split semisimple matrices. Hence $A\in M_{2,3}\setminus Q^{ss}_{2,3}$.
\end{example}

\begin{proof}[Proof of Theorem \ref{thm7}]
As each generalized Jordan block $J_{f,m}$ is similar to a matrix $\mathcal A$ as in Equation (\ref{EQ3}) but of size $s\times s$ with $s:=m\deg(f)$, it suffices to write $\mathcal A=B+C$ as a sum of two $k$th powers. Let $t:=\Tr(\mathcal A)$. By Corollary \ref{cyclicLU}, we have a decompostion using $(1,0,0,1,\ldots,1,0,0,1,c,t-s+1-c)\in\bbf_q^{2s}$ if $s$ is odd, and $(0,1,1,0,\ldots,0,1,1,0,0,1,c,t-s+1-c)\in\bbf_q^{2s}$ if $s$ is even, where $c\in\bbf_q$ is to be determined. If $t\neq s-1$, then $t-s+1\neq 0$, and it suffices to write $t-s+1\in\bbf_q^{\times}$ as a sum of two $d$-th powers $t-s+1=\beta^d+\gamma^d$ for some $\beta,\gamma\in\bbf_q$ such that $\beta\gamma\neq 0,\gamma^d\neq 1$. In other words, $c=\beta^d$. This can be done if $q\ge (d-1)^4+6d$ by Corollary \ref{Jolyc} applied to $(\eta_1,\eta_2)=(1,2)$. On the other hand, if $t=s-1$, then we change the $(s-1,s-1)$ entry of $B$ from $0$ to $1$, hence the sum of the $(s,s)$ entries of $B$ and $C$ is $-1\neq 0$. Now we need $-1=\beta^d+\gamma^d$ for some $\beta,\gamma\in\bbf_q$ such that $\beta^d\neq 1,\gamma^d\not\in\{0,1\}$. Again this can be done if $q\ge (d-1)^4+6d$.

So, as in the proof of Theorem \ref{thm2} (b), it follows that $P^{ss}=M_{n,q}$.
\end{proof}

\begin{proof}[Proof of Theorem \ref{thm8}]
Considering a cyclic decomposition of a matrix with entries in $\mathbb F_q$, it suffices to write each of the following kinds of matrices  as sums of two $k$th powers: an $n\times n$ nonscalar matrix with $n\ge 2$ and an $m\times m$ scalar matrix with $m\in\{2,3\}$.

Let $A\in M_{n,q}$ be a nonscalar matrix and let $t:=\Tr(A)-2(n-2)$. Suppose $t\neq 0$ and $q\ge \sqrt[3]{(d-1)^{8}+8d}$. Then 
$$N_4(t)\ge q^3-(d-1)^4q^{3/2}=\frac{q^{3/2}}{q^{3/2}+(d-1)^4}(q^3-(d-1)^8)>2^{-1}8d=4d$$ 
by Theorem \ref{Joly} applied to $s=4$ and $(k_i,a_i)=(d-1)$ for $i\in\{1,\ldots,4\}$, and hence we can write $t=a_1^d+a_2^d+a_3^d+a_4^d$ with $a_3a_4\neq 0$, as the number of solutions of the equation $t=x_1^d+x_2^d+x_3^d+x_4^d$ with three variables equal to $0$ is at most $4d$. Suppose $t=0$ and $q\ge \sqrt{(d-1)^6+6d}$. Then $N_3(t-1)\ge q^2-(d-1)^3q>3d$ by Theorem \ref{Joly} applied to $(s,k_1,\ldots,k_s,a_1,\ldots,a_s)=(3,d,d,d,1,1,1)$, and we can write $t=a_1^d+a_2^d+a_3^d+a_4^d$ with $a_1=1, a_3a_4\neq 0$, as the number of solutions of the equation $t-1=x_1^d+x_2^d+x_3^d$ with two variables equal to $0$ is at most $3d$. For $i\in\{1,2,3,4\}$ let $b_i\in\bbf_q$ be such that $a_i^d=b_i^k$. Based on Theorem \ref{nonscalar} and Corollary \ref{generalLU}, we have $A\in P_{k,n,q}$ by decomposing $A=L+U$ where $L$ is lower triangular and its diagonal entries are $1,\ldots,1,b_1^k,b_3^k$, and $U$ is an upper triangular and its diagonal entries are $1,\ldots,1,b_2^k,b_4^k$.

Let $A:=aI_m\in M_{m,q}$ be a scalar matrix with $m\in\{2,3\}$. Then $A\in P_{k,n,q}$ by Corollary \ref{const} and the inequality $q\ge \lfloor(d_m-1)^{4/m}\rfloor+1$.

Thus the theorem holds if $q$ is greater than its maximum. Notice that $\sqrt{(d-1)^6+6d}\ge \sqrt[3]{(d-1)^{8}+8d}$ for $d\ge 1$.\end{proof}

\begin{proof}[Proof of Corollary \ref{cor8}]
As $\max\{d,d_2,d_3\}\le k$ and $(k-1)^2+1\le\sqrt{(k-1)^6+6k}$ for $k\ge 1$, it suffices to compare $k^3-3k^2+3k$ with $\sqrt{(k-1)^6+6k}$ and check the case when $k\le 2$. If $k\ge 4$, then $k^3-3k^2+3k\ge \sqrt{(k-1)^6+6k}$. The case $k=3$ follows from the set identity $\{\mathbb F_q|q\ge k^3-3k^2+3k\}=\{\mathbb F_q|q\ge \sqrt{(k-1)^6+6k}\}$ as $q$s are prime powers. The case $k=1$ is trivial.

Suppose $k=2$. As $\sqrt{(k-1)^6+6k}=\sqrt{13}$, it suffices to check $q=3$. Let $A\in M_{n,3}$. Considering a Jordan block decomposition of $A$, if a block is a scalar matrix, then it is a sum of two squares as $0=0^2+0^2$, $1=1^2+0^2$, and $2=1^2+1^2$. If the block is nonscalar, then by Theorem \ref{nonscalar}, it is similar to a matrix with diagonal entries equal to $2,2,\ldots,2,a$, where $a\in\bbf_3$. Hence by Corollary \ref{generalLU}, $A\in P_{k,n,3}$.
\end{proof}

\begin{proof}[Proof of Theorem \ref{thm9}]
The set $\mathbb P_{d,q}\setminus\{0\}$ is a cyclic group $\langle a\rangle$ of order $o:=(q-1)/d\ge 2$; so $a\neq 1$. As $\cup_{s=1}^{\infty} s\mathbb{P}_{d,q}=\bbf_{p^{\ell}}$, we have $\mathbb F_{p^{\ell}}=\mathbb F_p[a]$; so $a$ has degree $\ell$ over $\bbf_p$ and $a-1,\ldots,a^{\ell}-1$ are linearly independent over $\bbf_p$. Thus the order of $a$ is at least $\ell+1$, i.e., the inequality $o\ge\ell+1$ holds, and $\{a-1,\ldots,a^{\ell}-1\}$ is an $\bbf_p$-basis of $\bbf_{p^{\ell}}$. Let $m:=n-1-\ell(p-1)\ge 0$. 

We first assume that $o\ge \ell+2$ and $n>\ell(p-1)$; so $m\ge 0$. As $o\ge \ell+2$, there exists $b\in\bbp_{d,q}\setminus\{0,1,a,\ldots,a^{\ell}\}$. For each $A\in M_{n,q}$, there exists a unique $(c_1,\ldots,c_{\ell})\in\{0,\ldots,p-1\}^{\ell}$ such that 
\[
\Tr(A)=\sum_{i=1}^{\ell}(p-1-c_i)a^i+b+(m+c_1+\cdots+c_{\ell})\cdot 1.
\]
If $A$ is cyclic, we can decompose $\mathcal A=B+C$ as a sum of two $k$th powers of split semisimple matrices over $\bbf_q$ by Corollary \ref{cyclicLU} using the following vector in $\bbf_q^{2n}$ whose repeated sequences have period $4$:
\[
\left(\underbrace{a,0,0,a,\ldots}_{2(p-1-c_1)},\ldots,\underbrace{\ldots,a^{\ell},0,0,a^{\ell}}_{2(p-1-c_{\ell})},\underbrace{\ldots,0,1,1,0,\ldots}_{2(m+c_1+\cdots+c_{\ell})},0,b\right),
\]
where the numbers under the braces mean the number of elements in the braces. Note that if $n$ is odd, then all diagonal entries of $B$ (resp. $C$) with odd positions are zero (resp. nonzero); if $n$ is even, then all diagonal entries of $B$ (resp. $C$) with even positions are zero (resp. nonzero). Thus $\mathcal A, A\in P$. E.g., if $(n,q)=(6,5)$ (so $\ell=l=1$) and $c_1=2$, then
\[
\mathcal A=\begin{pmatrix}
    a&&&&&\\    
    1&0&&&&\\
    &&1&&&\\    
    &&1&0&&\\    
    &&&&1&\\
    &&&&1&0
\end{pmatrix}+
\begin{pmatrix}
    0&&&&&*\\
    &a&&&&*\\    
    &1&0&&&*\\
    &&&1&&*\\
    &&&1&0&*\\
    &&&&&b
\end{pmatrix}.
\]

We now assume that $o=\ell+1$ and $n>(\ell+1)(p-1)$; so $m\ge p-1$. In this case we uniquely write
\[
\Tr(A)=\sum_{i=1}^{\ell}(p-1-c_i)a^{i+1}+1+(m+c_1+\cdots+c_{\ell})\cdot a.
\]
Let $\epsilon\in\{0,1\}$ be such that $-\epsilon+\sum_{i=1}^{\ell-1} (p-1-c_i)$ is even. If $A$ is cyclic, we can decompose $\mathcal A=B+C$ by Corollary \ref{cyclicLU} using the following vector in $\bbf_q^{2n}$ whose repeated sequences have period $4$:
\[
\left(\underbrace{a^2,0,0,a^2,\ldots}_{2(p-1-c_1)},\ldots,\underbrace{\ldots,a^{\ell},0,0,a^{\ell}}_{2(p-1-c_{\ell-1})},\underbrace{1,0,0,a,\ldots}_{4(p-1-c_{\ell})},\underbrace{\ldots,0,a,a,0}_{N},0,1\right)
\]
if $n+\epsilon$ is even, and
\[
\left(\underbrace{a^2,0,0,a^2,\ldots}_{2(p-1-c_1)},\ldots,\underbrace{\ldots,0,a^{\ell},a^{\ell},0}_{2(p-1-c_{\ell-1})},\underbrace{0,a,1,0,\ldots}_{4(p-1-c_{\ell})},\underbrace{\ldots,0,a,a,0}_{N},0,1\right)
\]
if $n+\epsilon$ is odd, where $N:=2(m+c_1+\cdots+c_{\ell-1}+2c_{\ell}-p+1)$. Then we can decompose $\mathcal A=B+C$ in such a way that the only diagonal entry of $C$ equal to $1$ is on the $(n,n)$ position, and as in the previous case we argue that $\mathcal A, A\in P$.\end{proof}

\begin{proof}[Proof of Theorem \ref{thm9b}]
For distinct $x_1^k,x_2^k\in\mathbb F_p^{\times}\setminus\{1\}\cap \mathbb P_{d,q}$, let $c\in\{0,\ldots,p-1\}$ be such that $\Tr(A)=(p-1-c)x_1^k+x_2^k+(n-p+c)\cdot 1$. The remaining part of the proof is similar to the proof of Theorem \ref{thm9} in the case $o\ge\ell+2$. 
\end{proof}

\begin{theorem}\label{thm10}
Suppose $q=p>d+1$ and $p$ does not divide $k$. Then $P=M_{n,q}$ if 
$n>\max\left\{\frac{p-1}{2},\frac{4\ln{(k-1)}}{\ln{p}}\right\}$, and $Q=M_{n,q}$ if 
$n>\max\left\{\frac{p-1}{2},\frac{\ln{\bigl((k-1)^4+4k\bigr)}}{\ln{p}}\right\}$ and $p\frac{p^n-1}{d_n}$ is even.
\end{theorem}

\begin{proof}
If $A$ is a nonscalar matrix, then $A$ is quasi-cyclic. As $(p-1)/d\in\mathbb N\setminus\{1\}$, there exists $a\in (\bbf_p^{\times}\setminus\{1\})\cap\mathbb P_{d,p}$. Let $b:=\Tr(A)-2(n-\lceil{(p+1)/2}\rceil)\in\mathbb F_p$. As $1-a\neq 0$, there exists $x\in\bbf_p$ such that $b=(1-a)x+2\lceil{(p+1)/2}\rceil a$; so $b=x+(2\lceil{(p+1)/2}\rceil-x)a$. Let $r\in\{0,1,\ldots,p-1\}$ be such that $r+p\mathbb Z=x\in\bbf_p$. Let $u_i:=a$ for $1\le i\le (2\lceil{(p+1)/2}\rceil-r)$. As we are assuming $p\le 2n$, let $u_j:=1$ for $(2\lceil{(p+1)/2}\rceil-r)<j\le 2n$. Then $\Tr(A)=u_1+\cdots+u_{2n}$. Hence $A\in Q\subset P$ by the proof of Corollary \ref{generalLU}. 

If $p^n>(k-1)^4$, we have $\{aI_n|a\in\mathbb F_p\}\subset P$ by Corollary \ref{const}. If $p^n>(k-1)^4+4k$ and $p(p^n-1)/d_n$ is even, then $\{aI_n|a\in\mathbb F_p\}\subset Q$ by Corollary \ref{constQ}.
\end{proof}

\begin{theorem}\label{thm16}
If $q>d^2$, $p\nmid k$, and $p\in\{2,3\}$, then $P=M_{n,q}$ for 
$n>\max\left\{3,\frac{4\ln{(k-1)}}{\ln{q}}\right\}.$
\end{theorem}

\begin{proof}
Let $A\in M_{n,q}$. Suppose $A$ is a nonscalar matrix. As $q>d^2$, we have $\ell=l$ (see \cite[Section 3]{Win98}) and from \cite[Corollary 1]{Cip09} it follows that $\gamma(k,q)=\gamma(d,k)\le 8$. As $2n\ge 8$, there exist $u_1,\ldots,u_{2n}\in\mathbb F_q$ such that $\Tr(A)=\sum_{i=1}^{2n} u_i^k$. As $p$ zero elements in the sequence $u_1,\ldots,u_{2n}$ can be replaced by $p$ one elements and as $p\in\{2,3\}$, we can assume that $\prod_{i=1}^{2n-2} u_i\neq 0$. Based on this, Theorem \ref{nonscalar}, and Corollary \ref{generalLU}, it follows that $A\in P$. If $A$ is a scalar matrix, then $A\in P$ by Corollary \ref{const} and the inequality $q^n>(k-1)^4$.
\end{proof}

\begin{proof}[Proof of Theorem \ref{thm18}]
Let $a\in \bigl(\bbf_p^{\times}\setminus\{1\}\bigr)\cap\mathbb P_{d,q}$. Let $r\in\{2,\ldots,p-1\}$ be such that $r+p\mathbb Z=a$. Then $0=a+(p-r)\cdot 1=p\cdot 1\in\mathbb F_p$ is a sum of either $(p-r+1)$ or $p$ elements in $\mathbb P_{d,q}\setminus\{0\}$. As $p$ and $p-r+1$ are coprime to each other, each integer $m\ge (p-r)(p-1)$ can be written as $m=s_{1}(p-r+1)+s_{2}p$ with $s_1,s_2\in\bbz_{\ge 0}$ and thus $0\in\mathbb F_p$ can be written as a sum of $m$ elements in $\mathbb P_{d,q}\setminus\{0\}$. Hence by a similar argument as in the proof of Theorem \ref{thm16}, it suffices to note that if $A\in M_{n,q}$ is nonscalar with $n>(p^2-3p+7)/2$, we have $2n\ge 8+(p-1)(p-2)\ge \gamma(d,p)+(p-1)(p-2)$ and we can write $\Tr(A)=\sum_{i=1}^s u_i^k=\sum_{i=1}^{2n} u_i^k$, where $s\in\{2,\ldots,8\}$ (so $2n-s\ge (p-2)(p-1)$), and $u_1,\ldots,u_{2n}\in\mathbb F_q$ with $u_i\neq 0$ if $i\notin\{s-1,s\}$.
\end{proof}

\begin{proof}[Proof of Theorem \ref{nlarge}]
Let $\gamma^*(d,q)\in\mathbb N$ be the minimum of the set $\{s\in\mathbb N|\bbf_q=s(\mathbb P_{k,d,q}\setminus\{0\})\}$. It exists as $\ell=l$ and from \cite[Theorem 1.1]{CC12} we have $\gamma^*(d,q)\le d+1$ (resp. $\gamma^*(d,q)=2d$). An easy induction shows that we have $\bbf_q=t(\mathbb P_{k,d,q}\setminus\{0\})$ for all integers $t\ge\gamma^*(d,q)$. Hence by a similar argument as in the proof of Theorem \ref{thm16}, it suffices to note that if $A\in M_{n,q}$ is a nonscalar matrix with $2n\ge d+1$ (resp. $n\ge d$), then we can write $\Tr(A)=\sum_{i=1}^{2n} u_i^k$ with $u_1,\ldots,u_{2n}\in\mathbb F_q^{\times}$.
\end{proof}

\begin{example}\label{EX4}
Assume $k=2$. We check that $P=M_{n,q}$ if $q>2$. If $p=2$, then $d=1$ and from Theorem \ref{thm2} (b) it follows that for $q>2$ we have $Q^{ss}=M_{n,q}$. Thus we can assume that $p>2$, hence $p$ does not divide $k$, and the equality $P=M_{n,q}$ follows from Corollary \ref{cor8}. 
\end{example}

\begin{example}\label{EX5}
Assume $k=3$. We check that $P=M_{n,q}$ if $q\notin\{4,7\}$. If $p=3$, then $d=1$ and from Theorem \ref{thm2} (b) it follows that $P=M_{n,q}$. Thus we can assume that $p\neq 3$, i.e., $p$ does not divide $k$, and hence for $q\notin\{4,7\}$ the equality $P=M_{n,q}$ follows from Corollary \ref{cor8} if $q\ge 11$, from Theorem \ref{thm3} if $q=2$, and from Theorem \ref{thm2} (b) if $q\in\{5,8\}$. 

We assume $q=4$. Thus $(\ell,l)=(1,2)$ and we have $d_m=3$ for all $m\in\mathbb N$. From Example \ref{EX3} if follows that all nonscalar matrices in $M_{n,q}$ of trace in $\mathbb F_2$ belong to $P$. From Corollary \ref{const} applied with $s=1$ it follows that for $n\ge 3$ the scalar matrices in $M_{n,4}$ belong to $P$.

We assume $q=7$. Thus $\ell=l=1$ and we have $d_m=3$ for all $m\in\mathbb N$. From Theorem \ref{nlarge} it follows that $P=M_{n,q}$ for $n\ge 6$. We assume now that $n\in\{2,3,4,5\}$. From Corollary \ref{const} applied with $s=1$ it follows that the scalar matrices in $M_{n,7}$ belong to $P$. From Example \ref{EX3} if follows that all nonscalar matrices of trace in the set $\{2n-2,2n-1,2n\}$ belong to $P$.
\end{example}

\section{The case of \texorpdfstring{$\mathbb F_2$}{F2}}\label{S6}

In this section we assume that $q=2$. We first treat the case $n=2$.

The 16 elements of $M_{2,2}$ are grouped in conjugacy classes as follows:

\noindent
(i) One zero matrix: $0_2$.

\noindent
(ii) Three nonzero nilpotent matrices: $\begin{pmatrix}
    0 & 1\\ 0 & 0
\end{pmatrix}, \begin{pmatrix}
    0 & 0\\ 1 & 0
\end{pmatrix},\begin{pmatrix}
    1 & 1\\ 1 & 1
\end{pmatrix}$.

\noindent
(iii) One identity element: $I_2$.

\noindent
(iv) Six idempotents of rank 1: $\begin{pmatrix}
    1 & 0\\ 0 & 0
\end{pmatrix},\begin{pmatrix}
    1 & 0\\ 1 & 0
\end{pmatrix},\begin{pmatrix}
    1 & 1\\ 0 & 0
\end{pmatrix}$, and their additive translates by $I_2$.

\noindent
(v) Two invertible matrices of order 3: $\begin{pmatrix}
    1 & 1\\ 1 & 0
\end{pmatrix}, \begin{pmatrix}
    0 & 1\\ 1 & 1
\end{pmatrix}$.

\noindent
(vi) Three invertible matrices of order 2: $\begin{pmatrix}
    1 & 1\\ 0 & 1
\end{pmatrix}, \begin{pmatrix}
    1 & 0\\ 1 & 1
\end{pmatrix},\begin{pmatrix}
    0 & 1\\ 1 & 0
\end{pmatrix}$.

\begin{theorem}\label{M22}
Let $k\ge 2$. Then
\[
P=\begin{cases}
M_{2,2} & {\rm\ if}\ k\not\equiv 0\ {\rm mod}\ 6\\
M_{2,2}\setminus\{{\rm the\ two\ invertible\ matrices\ of\ order\ 3}\} & {\rm\ if}\ k\equiv 0\ {\rm mod}\ 6.
\end{cases}
\]

\end{theorem}

\begin{lemma}\label{idemp22}
We have $\Pi_{2,2}=M_{2,2}\setminus\{{\rm the\ invertible\ matrices\ of\ order\ 3}\}$.
\end{lemma}
\noindent{\it Proof.}
As $\mathbb{I}_{2,2}$ is stable under conjugation, we only need to check this for one representative of a conjugacy class. We have the following identities $0_2=0_2+0_2=I_2+I_2$, $\begin{pmatrix}
    0 & 1\\ 0 & 0
\end{pmatrix}=\begin{pmatrix}
    1 & 0\\ 0 & 0
\end{pmatrix}+\begin{pmatrix}
    1 & 1\\ 0 & 0
\end{pmatrix}$, $I_2=0_2+I_2$, $\begin{pmatrix}
    1 & 0\\ 0 & 0
\end{pmatrix}=0_2+\begin{pmatrix}
    1 & 0\\ 0 & 0
\end{pmatrix}$, $\begin{pmatrix}
    1 & 1\\ 0 & 1
\end{pmatrix}=\begin{pmatrix}
    1 & 1\\ 0 & 0
\end{pmatrix}+\begin{pmatrix}
    0 & 0\\ 0 & 1
\end{pmatrix}$.

As $\begin{pmatrix}
    1 & 1\\ 1 & 0
\end{pmatrix}$ has trace $1$, if it is a sum of two idempotents, one of them has to have trace $0$. Hence it must be $0_2$ or $I_2$, which is not possible as $\begin{pmatrix}
    1 & 1\\ 1 & 0
\end{pmatrix}$ and $\begin{pmatrix}
    1 & 1\\ 1 & 0
\end{pmatrix}-I_2$ do not belong to $\mathbb{I}_{2,2}$. $\square$

\vskip 2mm

\noindent{\it Proof of Theorem \ref{M22}.}
As $\Pi_{2,2}\subset P_{k,2,2}$, from Lemma \ref{idemp22} it follows that 
\[
M_{2,2}\setminus\{{\rm the\ two\ invertible\ matrices\ of\ order\ 3}\}\subset P.
\]
To test if invertible matrices of order $3$ do or do not belong to $P$, as $|G_{2,2}|=6$, by working modulo $6$ we can assume that $1\le k\le 6$. By Theorem \ref{thm3}, if $k$ is odd, then $P=M_{2,2}$. If $k\in\{2,4\}$, then $(v)\subset\mathbb P_{k,2,2}$. If $k=6$, then $\mathbb P_{k,2,2}=\mathbb{I}_{2,2}$ and $P_{k,2,2}=\Pi_{2,2}\subsetneq M_{2,2}$.\ \ $\square$

\begin{remark}\label{remark4}
Lemma \ref{idemp22} implies that $L_{\Pi_{2,2}}=\{I_2\}$.

If $k=3$, then $G_{2,2}\cap\mathbb P_{3,2,2}=(iii)\cup (vi)$, hence $I_2\not\in Q$ and thus Theorem \ref{thm2} is not true for $q=2$.

As each nonzero matrix $A\in M_{2,2}$ is uniquely (up to order) a sum of two invertible matrices (see \cite[Theorem]{Lor87}), by taking $A=I_2+\begin{pmatrix}
    1 & 1\\ 0 & 1\end{pmatrix}$, as $I_2$ is not cyclic and $\begin{pmatrix}
    1 & 1\\ 0 & 1\end{pmatrix}$ is not semisimple, it follows that Theorem \ref{thm5} (a) and (b) are not true for $q=2$. 

As $\Pi_{2,2}=P^{ss}_{k,2,2}$ and the non-split semisimple elements are in $(v)$ and their $k$th powers are in $(v)$ (resp. in $(iii)$) if $\gcd(k,3)=1$ (resp. if $3|k$), it follows that $P^{s}_{k,2,2}:=\{A^k+B^k\mid {\rm semisimple}\ A,B\in M_{2,2}\}$ equals $M_{2,2}$ (resp. $\Pi_{2,2}$).
\end{remark}

Next we consider the case $n=3$.

The 512 elements of $M_{3,2}$ are grouped in conjugacy classes as follows:

\noindent
(1) One zero matrix: $M_1:=0_3$.

\noindent
(2) One identity matrix: $M_2:=I_3$.

\noindent
(3) Twenty-eight idempotent matrices conjugate to $M_3:=\begin{pmatrix}
    1 & 0 & 0\\ 0 & 0 & 0\\ 0& 0& 0
\end{pmatrix}$.

\noindent
(4) Twenty-eight idempotent matrices conjugate to $M_4:=\begin{pmatrix}
    1 & 0 & 0\\ 0 & 1 & 0\\ 0& 0& 0
\end{pmatrix}$.

\noindent
(5) Eighty-four matrices conjugate to $M_5:=\begin{pmatrix}
    1 & 0 & 0\\ 0 & 0 & 1\\ 0 & 0 & 0
\end{pmatrix}$; so $M_5^2=M_3$.

\noindent
(6) Twenty-one nilpotent matrices conjugate to $M_6:=\begin{pmatrix}
    0 & 1 & 0\\ 0 & 0 & 0\\ 0& 0& 0
\end{pmatrix}$.

\noindent
(7) Forty-two nilpotent matrices conjugate to $M_7:=\begin{pmatrix}
    0 & 1 & 0\\ 0 & 0 & 1\\ 0& 0& 0
\end{pmatrix}$.

\noindent
(8) Eighty-four matrices conjugate to $M_8:=\begin{pmatrix}
    1 & 1 & 0\\ 0 & 1 & 0\\ 0 & 0 & 0
\end{pmatrix}$; so $M_8^2=M_4$.

\noindent
(9) Twenty-one matrices of order $2$ conjugate to $M_9:=\begin{pmatrix}
    1 & 1 & 0\\ 0 & 1 & 0\\ 0& 0& 1
\end{pmatrix}$.

\noindent
(10) Forty-two matrices of order $4$ conjugate to $M_{10}:=\begin{pmatrix}
    1 & 1 & 0\\ 0 & 1 & 1\\ 0& 0& 1
\end{pmatrix}$.

\noindent
(11) Fifty-six matrices conjugate to $M_{11}:=\begin{pmatrix}
    0 & 1 & 0\\ 1 & 1 & 0\\ 0 & 0 & 0
\end{pmatrix}$; so $M_{11}^3=M_4$.

\noindent
(12) Fifty-six matrices of order $3$ conjugate to $M_{12}:=\begin{pmatrix}
    0 & 1 & 0\\ 1 & 1 & 0\\ 0& 0& 1
\end{pmatrix}$.

\noindent
(13) Twenty-four matrices of order $7$ conjugate to $M_{13}:=\begin{pmatrix}
    0 & 0 & 1\\ 1 & 0 & 1\\ 0& 1& 0
\end{pmatrix}$.

\noindent
(14) Twenty-four matrices of order $7$ conjugate to $M_{14}:=\begin{pmatrix}
    0 & 0 & 1\\ 1 & 0 & 0\\ 0& 1& 1
\end{pmatrix}$.

\begin{theorem}\label{M32}
Let $k\ge 2$. Then
\[
P=\begin{cases}
M_{3,2} & {\rm\ if\ } k\not\equiv 0\ {\rm mod}\ 42\\
M_{3,2}\setminus(14) & {\rm\ if\ } k\equiv 42\ {\rm mod}\ 84\\
M_{3,2}\setminus\bigl({ (11)\cup (12)\cup (13)\cup (14)}\bigr) & {\rm\ if\ } k\equiv 0\ {\rm mod}\ 84.
\end{cases}
\]
\end{theorem}

\begin{lemma}\label{idemp32}
We have $\Pi_{3,2}=M_{3,2}\setminus\bigl({ (11)\cup (12)\cup (13)\cup(14)}\bigl)$.
\end{lemma}

\begin{proof}
Clearly, the idempotent matrices (1), (2), (3), and (4) belong to $\Pi_{3,2}$. We have the following identities:

$M_5=\begin{pmatrix}
    0 & 0 & 0\\ 0 & 0 & 0\\ 0& 0& 1
\end{pmatrix}+\begin{pmatrix}
    1 & 0 & 0\\ 0 & 0 & 1\\ 0& 0& 1
\end{pmatrix}$, 
$M_6=\begin{pmatrix}
    1 & 0 & 0\\ 0 & 0 & 0\\ 0& 0& 0
\end{pmatrix}+\begin{pmatrix}
    1 & 1 & 0\\ 0 & 0 & 0\\ 0& 0& 0
\end{pmatrix}$,

$M_7=\begin{pmatrix}
    0 & 1 & 0\\ 0 & 1 & 0\\ 0& 0& 0
\end{pmatrix}+\begin{pmatrix}
    0 & 0 & 0\\ 0 & 1 & 1\\ 0& 0& 0
\end{pmatrix}$, 
$M_8=\begin{pmatrix}
    1 & 0 & 0\\ 0 & 0 & 0\\ 0& 0& 0
\end{pmatrix}+\begin{pmatrix}
    0 & 1 & 0\\ 0 & 1 & 0\\ 0& 0& 0
\end{pmatrix}$,

$M_9=\begin{pmatrix}
    1 & 0 & 0\\ 0 & 0 & 0\\ 0& 0& 0
\end{pmatrix}+\begin{pmatrix}
    0 & 1 & 0\\ 0 & 1 & 0\\ 0& 0& 1
\end{pmatrix}$, 
$M_{10}=\begin{pmatrix}
    0 & 1 & 0\\ 0 & 1 & 0\\ 0& 0& 0
\end{pmatrix}+\begin{pmatrix}
    1 & 0 & 0\\ 0 & 0 & 1\\ 0& 0& 1
\end{pmatrix}$.

Let $j\in\{11,12,13,14\}$. As $\mathbb I_{3,2}=\cup_{i=1}^ 4 (i)$ and as $(2)\cup (4)$ is just the translation of $(1)\cup (3)$ by $I_3$, to show that $M_j\notin\Pi_{3,2}$, it suffices to show that $M_j$ and $M_j+I_3$ cannot be written as the sum of two elements  in $(1)\cup (3)$, which can be easily checked. Hence $M_j\notin\Pi_{3,2}$.
\end{proof}

\begin{lemma}\label{32class3}
If $42\nmid k$, then $P=M_{3,2}$.
\end{lemma}

\begin{proof}
 If $2\nmid k$, then $P=M_{3,2}$ by Theorem \ref{thm3}. Hence we assume $2|k$.
 
 If $3\nmid k$, then $A:=\begin{pmatrix}
    0 & 1 \\ 1 & 1 
\end{pmatrix}\in P_{k,2,2}$ by Theorem \ref{M22}, hence $M_{11}=$\\[-1pt]$\textup{Diag}(A,0),M_{12}=\textup{Diag}(A,1)\in P_{k,3,2}$; moreover, $M_{13}=\begin{pmatrix}
    1 & 0 & 0\\ 0 & 0 & 0\\ 0 & 0 & 0
\end{pmatrix}+$\\[-1pt]$\begin{pmatrix}
    1 & 0 & 1\\ 1 & 0 & 1\\ 0 & 1 & 0
\end{pmatrix}$ and $M_{14}=\begin{pmatrix}
    0 & 0 & 0\\ 0 & 0 & 0\\ 0 & 0 & 1
\end{pmatrix}+\begin{pmatrix}
    0 & 0 & 1\\ 1 & 0 & 0\\ 0 & 1 & 0
\end{pmatrix}$, where $\begin{pmatrix}
    1 & 0 & 1\\ 1 & 0 & 1\\ 0 & 1 & 0
\end{pmatrix}$ and $\begin{pmatrix}
    0 & 0 & 1\\ 1 & 0 & 0\\ 0 & 1 & 0
\end{pmatrix}$ are $k$th powers of elements in (11) and (12) respectively, and $\begin{pmatrix}
    1 & 0 & 0\\ 0 & 0 & 0\\ 0 & 0 & 0
\end{pmatrix},\begin{pmatrix}
    0 & 0 & 0\\ 0 & 0 & 0\\ 0 & 0 & 1
\end{pmatrix}\in\mathbb{I}_{3,2}$.

If $7\nmid k$, then $M_{11}=\begin{pmatrix}
    0 & 1 & 1\\ 0 & 0 & 0\\ 0 & 1 & 1
\end{pmatrix}+\begin{pmatrix}
    0 & 0 & 1\\ 1 & 1 & 0\\ 0 & 1 & 1
\end{pmatrix}$ and $M_{12}=\begin{pmatrix}
    0 & 1 & 1\\ 0 & 0 & 0\\ 0 & 1 & 1
\end{pmatrix}+$\\[-1pt]$\begin{pmatrix}
    0 & 0 & 1\\ 1 & 1 & 0\\ 0 & 1 & 0
\end{pmatrix}$, where $\begin{pmatrix}
    0 & 0 & 1\\ 1 & 1 & 0\\ 0 & 1 & 1
\end{pmatrix}$ and $\begin{pmatrix}
    0 & 0 & 1\\ 1 & 1 & 0\\ 0 & 1 & 0
\end{pmatrix}$ are $k$th powers of elements in (13) or (14) as they have order $7$ and $\begin{pmatrix}
    0 & 1 & 1\\ 0 & 0 & 0\\ 0 & 1 & 1
\end{pmatrix}\in\mathbb{I}_{3,2}$; moreover, we have $M_{13},M_{14}\in P$ as they have order $7$.
\end{proof}

\begin{lemma}\label{32class2}
If $42|k$, then  $(14)\cap P=\emptyset$.
\end{lemma}

\begin{proof}
We have $\mathbb P_{k,3,2}=\mathbb{I}_{3,2}\cup(9)$ and $\Tr(M_9)=\Tr(M_{14})=1$, hence by reasons of traces and Lemma \ref{idemp32} it follows that $M_{14}\notin P$.
\end{proof}

\begin{lemma}\label{32class1}
Suppose $k=42m$ with $m\in\mathbb{N}$. For $j\in\{11,12,13\}$, the conjugacy class ($j$) is contained in $P$ if and only if $m$ is odd.
\end{lemma}
\begin{proof}
   
If $m$ is odd, then $k\equiv 6$ mod $12$. We have $M_{11}=\begin{pmatrix}
    1 & 1 & 0\\ 0 & 1 & 0\\ 0 & 0 & 1
\end{pmatrix}+\begin{pmatrix}
    1 & 0 & 0\\ 1 & 0 & 0\\ 0 & 0 & 1
\end{pmatrix}$, $M_{12}=\begin{pmatrix}
    1 & 1 & 0\\ 0 & 1 & 0\\ 0 & 0 & 1
\end{pmatrix}+\begin{pmatrix}
    1 & 0 & 0\\ 1 & 0 & 0\\ 0 & 0 & 0
\end{pmatrix}$, and $M_{13}=\begin{pmatrix}
    0 & 0 & 0\\ 0 & 1 & 0\\ 0 & 0 & 0
\end{pmatrix}+\begin{pmatrix}
    0 & 0 & 1\\ 1 & 1 & 1\\ 0 & 1 & 0
\end{pmatrix}$, where $\begin{pmatrix}
    1 & 1 & 0\\ 0 & 1 & 0\\ 0 & 0 & 1
\end{pmatrix}$, $\begin{pmatrix}
    0 & 0 & 1\\ 1 & 1 & 1\\ 0 & 1 & 0
\end{pmatrix}$ are $k$th powers of elements in (10), and $\begin{pmatrix}
    1 & 0 & 0\\ 1 & 0 & 0\\ 0 & 0 & 1
\end{pmatrix}, \begin{pmatrix}
    1 & 0 & 0\\ 1 & 0 & 0\\ 0 & 0 & 0
\end{pmatrix}$,$\begin{pmatrix}
    0 & 0 & 0\\ 0 & 1 & 0\\ 0 & 0 & 0
\end{pmatrix}\in\mathbb{I}_{3,2}$.

If $2|m$, then $\mathbb P_{k,3,2}=\mathbb{I}_{3,2}$ by Proposition \ref{exponent}; so $M_{11},M_{12},M_{13}\notin\Pi_{3,2}$ by Lemma \ref{idemp32}.\end{proof}

\begin{proof}[Proof of Theorem \ref{M32}]
It follows from Lemmas \ref{idemp32}, \ref{32class3}, \ref{32class2}, and \ref{32class1}.
\end{proof}

\section{The case of \texorpdfstring{$\mathbb F_3$}{F3}}\label{S7}

In this section we assume $q=3$.

\begin{prop}\label{prop75}
If $6$ does not divide $k$, then $P=M_{n,3}$.
\end{prop}

\begin{proof}
If $d=1$, then this follows from Theorem \ref{thm2} (b). If $d=2$, then $3$ does not divide $k$. As $\bbf_3=2\mathbb{P}_{k,3}$, scalar matrices in $M_{n,3}$ belong to $P$. If $A\in M_{n,3}$ is nonscalar, then $A\in P$ by Example \ref{EX3}. 
\end{proof}

The $81$ elements of $M_{2,3}$ are grouped in conjugacy classes as follows:

\noindent
\tcircle{1} One zero matrix: $N_1:=0_2$.

\noindent
\tcircle{2} One identity matrix: $N_2:=I_2$.

\noindent
\tcircle{3} Twelve idempotent matrices conjugate to $N_3:=\begin{pmatrix}
    0 & 0 \\ 0 & 1 
\end{pmatrix}$.

\noindent
\tcircle{4} Eight matrices of order $3$ conjugate to $N_4:=\begin{pmatrix}
    1 & 1 \\ 0 & 1 
\end{pmatrix}$.

\noindent
\tcircle{5} Eight nilpotent matrices conjugate to $N_5:=\begin{pmatrix}
    0 & 1 \\ 0 & 0 
\end{pmatrix}$.

\noindent
\tcircle{6} Twelve matrices conjugate to $N_6:=2N_3$; so $N_6^2=N_3$.

\noindent
\tcircle{7} Twelve matrices of order $2$ conjugate to $N_7:=\begin{pmatrix}
    1 & 0 \\ 0 & 2 
\end{pmatrix}$.

\noindent
\tcircle{8} One matrix of order $2$ conjugate to $N_8:=2I_2$.

\noindent
\tcircle{9} Eight matrices of order $6$ conjugate to $N_9:=\begin{pmatrix}
    2 & 1 \\ 0 & 2 
\end{pmatrix}$.

\noindent
\tcircle{10} Six matrices of order $4$ conjugate to $N_{10}:=\begin{pmatrix}
    0 & 2 \\ 1 & 0 
\end{pmatrix}$; so $N_{10}^2=2I_2$.

\noindent
\tcircle{11} Six matrices of order $8$ conjugate to $N_{11}:=\begin{pmatrix}
    0 & 1 \\ 1 & 2 
\end{pmatrix}$; so $N_{11}^4=2I_2$.

\noindent
\tcircle{12} Six matrices of order $8$ conjugate to $N_{12}:=\begin{pmatrix}
    0 & 1 \\ 1 & 1 
\end{pmatrix}$; so $N_{12}^4=2I_2$.

\begin{theorem}\label{M23}
Let $k\ge 2$. Then
\[
P=\begin{cases}
M_{2,3} & {\rm\ if\ } k\not\equiv 0,6,12,18\ {\rm mod}\ 24\\
M_{2,3}\setminus  \textnormal{\tcircle{5}} & {\rm\ if\ } k\equiv 6,18\ {\rm mod}\ 24\\
M_{2,3}\setminus \bigl(\textnormal{\tcircle{5}}\cup\textnormal{ \tcircle{9}}\cup \textnormal{\tcircle{10}}\cup \textnormal{\tcircle{12}}\bigr) & {\rm\ if\ } k\equiv 0\ {\rm mod}\ 12.
\end{cases}
\]

\end{theorem}

\begin{lemma}\label{idemp23}
We have $\Pi_{2,3}=M_{2,3}\setminus\bigl(\textnormal{\tcircle{5}}\cup\textnormal{ \tcircle{9}}\cup \textnormal{\tcircle{10}}\cup \textnormal{\tcircle{12}}\bigr)$.
\end{lemma}

\begin{proof}
Clearly, the idempotent matrices from \tcircle{1}, \tcircle{2}, and \tcircle{3} belong to $\Pi_{2,3}$. We have identities $N_4=\begin{pmatrix}
    1 & 0 \\ 0 & 0 
\end{pmatrix}+\begin{pmatrix}
    0 & 1 \\ 0 & 1 
\end{pmatrix}$,
$N_6=N_3+N_3$, $N_7=\begin{pmatrix}
    0 & 0 \\ 0 & 1 
\end{pmatrix}+I_2$,
$N_8=I_2+I_2$, $N_{11}=\begin{pmatrix}
    0 & 1 \\ 0 & 1 
\end{pmatrix}+\begin{pmatrix}
    0 & 0 \\ 1 & 1 
\end{pmatrix}$. For $j\in\{5,9,10,12\}$, as $\Tr(N_j)\in\{0,1\}$, it is easy to see that $N_j\notin\Pi_{2,3}$.\end{proof}

\begin{lemma}\label{23class1}
If $6$ divides $k$, then $\textnormal{\tcircle{5}}\cap P=\emptyset$.
\end{lemma}
\begin{proof}
We have $\mathbb P_{k,2,3}=\mathbb{I}_{2,3}\cup\tcircle{8}\cup\tcircle{10}$. As $\Tr(N_5)=\Tr(N_{10})=0$ and $\Tr(N_8)=1$, by reasons of traces the assumption that $N_5=B+C$ with $B,C\in\mathbb P_{k,2,3}$ implies that the set $\{B,C\}\cap\{0_2,I_2\}$ is nonemepty and this leads to a contradiction; hence $N_5\notin P$.
\end{proof}

\begin{lemma}\label{23class2}
Suppose $k=6m$ with $m\in\mathbb N$. For $j\in\{9,10,12\}$, the conjugacy class \textnormal{\tcircle{j}} is contained in $P$ if and only if $m$ is odd.
\end{lemma}
\begin{proof}
If $m$ is odd, then we decompose $N_9=\begin{pmatrix}
    0 & 0 \\ 2 & 1 
\end{pmatrix}+\begin{pmatrix}
    2 & 1 \\ 1 & 1 
\end{pmatrix}$ and $N_{12}=\begin{pmatrix}
    1 & 0 \\ 0 & 0 
\end{pmatrix}+\begin{pmatrix}
    2 & 1 \\ 1 & 1 
\end{pmatrix}$, where $\begin{pmatrix}
    2 & 1 \\ 1 & 1
\end{pmatrix}\in\tcircle{10}$ and $N_{10}$ are $k$th powers of elements in $\tcircle{11}\cup\tcircle{12}$ and $\begin{pmatrix}
    0 & 0 \\ 2 & 1 
\end{pmatrix}, \begin{pmatrix}
    1 & 0 \\ 0 & 0 
\end{pmatrix}\in\mathbb{I}_{2,3}$.

If $m$ is even, then $\mathbb P_{k,2,3}=\mathbb{I}_{2,3}\cup\tcircle{8}$, hence $P_{k,3,2}=\Pi_{2,3}\cup (\mathbb{I}_{2,3}+2I_2)$. Clearly $(N_j-2I_2)\notin\mathbb I_{2,3}$. From the last two sentences and Lemma \ref{idemp23} it follows that $N_j\notin P_{k,3,2}$.\end{proof}

\begin{proof}[Proof of Theorem \ref{M23}]
Based on Proposition \ref{prop75}, to prove Theorem \ref{M23} we can assume that $k=6m$ with $m\in\mathbb N$. In this case, Theorem \ref{M23} follows from Lemmas \ref{idemp23}, \ref{23class1}, and  \ref{23class2}.
\end{proof}

\section{Open Problems}\label{S8}

\noindent
{\bf Problem 1 (Larsen's Conjecture).} Prove that for a given $k$, there exists a constant $N(k)\in\mathbb N$ depending only on $k$, such that for all pairs $(n,q)$, if
$q^{n^2}\ge N(k)$, then $P=M_{n,q}$.

\smallskip\noindent
{\bf Problem 2.} Classify all pairs $(n,q)$ for which there exists $d_{n,q}\in\mathbb N$ such that for  $k\in\mathbb N$ we have $P_{k,n,q}\neq M_{n,q}$ if and only if $d_{n,q}$ divides $k$.

\smallskip\noindent
{\bf Problem 3.} Classify all triples $(n,q,s)$ with $s$ a prime divisor of $e_{n,q}$ such that $\{A+B|A^{1+s^m}=A,B^{1+s^t}=B\; \textup{for some}\; t,m\in\mathbb N\}=M_{n,q}$.

\smallskip\noindent
{\bf Problem 4.} Classify all triples $(k,n,q)$ with $k<n$ and $P=2\mathbb P^{-}_{k,n,q}$ (see Remark \ref{remark1}).

\smallskip\noindent
{\bf Problem 5.} Generalize Lemma \ref{q3} and Theorem \ref{nonscalar} to larger classes of rings that contain the fields.

\smallskip\noindent
{\bf Problem 6.} If $P_{q,n,k}=M_{n,q}$, then also $P_{q,n,d}=M_{n,q}$. Determine when the converse holds.

\smallskip\noindent
{\bf Problem 7.} Classify all pairs $(k,q)$ such that $p\nmid k$ and the set
$N_{k,q}:=\{m\in\mathbb N|m\ge 2, Q_{k,m,q}=M_{m,q}\}$ is stable under addition.

\vskip 2mm
\hbox{Krishna Kishore}

\hbox{Department of Mathematics,} 

\hbox{Indian Institute of Technology (IIT) Tirupati,}

\hbox{Tirupati-Renigunta Road, Settipalli Post, Tirupati,}

\hbox{Andhra Pradesh, India 517506.}

\hbox{Email: kishorekrishna@iittp.ac.in}

\vskip 2mm
\hbox{Adrian Vasiu}

\hbox{Department of Mathematics and Statistics, Binghamton University,}

\hbox{P. O. Box 6000, Binghamton, New York 13902-6000, U.S.A.}

\hbox{Email: avasiu@binghamton.edu}

\vskip 2mm
\hbox{Sailun Zhan}

\hbox{Department of Mathematics and Statistics, Binghamton University,}

\hbox{P. O. Box 6000, Binghamton, New York 13902-6000, U.S.A.}

\hbox{Email: sailunzhan@gmail.com}


\begin{thebibliography}{DDDDD2}

\bibitem[Bot00]{Bot00} Botha, J. D.: Sums of diagonalizable matrices, {\it Linear Algebra Appl.} {\bf 315} (2000), no. 1-3, 1--23.

\bibitem[Bre18]{Bre18} Breaz, S.: Matrices over finite fields as sums of periodic and nilpotent elements, {\it Linear Algebra Appl.} {\bf 555} (2018), 92--97.

\bibitem[BCDM13]{BCDM13} Breaz, S.; Călugăreanu, G.; Danchev, P.; Micu, T.: Nil-clean matrix rings, {\it Linear Algebra Appl.} {\bf 439} (2013), no. 10, 3115--3119.

\bibitem[BM16]{BM16} Breaz, S.; Modoi, G. C.: Nil-clean companion matrices, {\it Linear Algebra Appl.} {\bf 489} (2016), 50--60.

\bibitem[Cip09]{Cip09} Cipra, J. A.: Waring's number in a finite field, {\it Integers} {\bf 9} (2009), A34, 435--440.

\bibitem[CC12]{CC12} Cochrane, T.; Cipra, J. A.: Sum-product estimates applied to Waring's problem over finite fields, {\it Integers} {\bf 12} (2012), no. 3, 385--403.

\bibitem[DeO73]{DeO73} De Oliveira, G. N.: Matrices with prescribed
entries and eigenvalues I, {\it Proc. Amer. Math. Soc.} {\bf 37}
(1973), no. 2, 380--386.

\bibitem[DK19]{DK19} Demiro$\check{g}$lu Karabulut, Y.: Waring's problem in finite
rings, {\it J. Pure Appl. Algebra} {\bf 223} (2019), no. 8, 3318--3329.

\bibitem[Fil69]{Fil69}
Fillmore, P. A.: On similarity and the diagonal of a matrix, {\it Amer. Math. Monthly} {\bf 76} (1969), 167--169.

\bibitem[FL58]{FL58} Farahat, H. K.; Ledermann, W.: Matrices with prescribed
characteristic polynomials, {\it Proc. Edinburgh Math. Soc.} {\bf 11},
1958/1959, 143--146.

\bibitem[GLS94]{GLS94} Gorenstein, D.; Lyons, R.; Solomon, R.: The
classification of the finite simple groups, Number 3, {\it Math. Surv.
and Monog.} {\bf 40}, American Mathematical Society, 1994.

\bibitem[Hen74]{Hen74} Henriksen, M.: Two classes of rings generated by their units, {\it J. Algebra} {\bf 31} (1974), 182--193.

\bibitem[Her83]{Her83} Hershkowitz, D.: Existence of matrices with prescribed
eigenvalues and entries, {\it Linear Multilinear Algebra} {\bf 14}
(1983), no. 4, 315--342.

\bibitem[HT88]{HT88} Hirano, Y.; Tominaga, H.: Rings in which every element is the sum of two idempotents, {\it Bull. Austral. Math. Soc.} {\bf 37} (1988), no. 2, 161--164.

\bibitem[Jol73]{Jol73} Joly, J. R.: Equations et vari\'et\'es alg\'ebriques sur un corps fini, {\it Enseignement Math.} {\bf 19} (1973), 1--117.

\bibitem[Kis22]{Kis22} Kishore, K.: Matrix Waring problem, {\it Linear Algebra
Appl.} {\bf 646} (2022), 84--94.

\bibitem[KS23]{KS23} Kishore, K.; Singh, A.: Matrix Waring problem II, https://arxiv.org

/abs/2205.04710. To appear in {\it Israel J. Math.}

\bibitem[Lan02]{Lan02} Lang, S.: Algebra, Rev. 3rd ed., {\it Grad. Texts Math.} {\bf 211}, Springer-Verlag New York, 2002.

\bibitem[Lor87]{Lor87} Lord, N. J.: Matrices as sums of invertible matrices,
{\it Math. Mag.} {\bf 60} (1987), no. 1, 33--35.

\bibitem[Mor06]{Mor06} Morrison, K. E.: Integer sequences and matrices
over finite fields, {\it J. Integer Seq.} Vol. {\bf 9} (2006), Article
06.2.1.

\bibitem[NP00]{NP00} Neumann, P. M.; Praeger, C. E.: Cyclic matrices in
classical groups over finite fields, Special issue in honor of Helmut
Wielandt, {\it J. Algebra} {\bf 234} (2000), no. 2, 367--418.

\bibitem[NVZ04]{NVZ04} Nicholson, W. K.; Varadarajan, K.; Zhou, Y.: Clean endomorphism rings, {\it Arch. Math. (Basel)} {\bf 83} (2004), no. 4, 340--343.

\bibitem[Pum07]{Pum07} Pumplün, S.: Sums of dth powers in non-commutative rings, {\it Beitr\"age Algebra Geom.} {\bf 48} (2007), no. 1, 291--301.

\bibitem[Rot95]{Rot95} Rotman, J. J.: An introduction to the theory of groups.
Fourth edition, {\it Grad. Texts Math.} {\bf 148}, Springer-Verlag, New
York, 1995.

\bibitem[Sma77]{Sma77} Small, C.: Sums of powers in large finite fields, {\it Proc. Amer. Math. Soc.} {\bf 65} (1977), 35--36.

\bibitem[Tor38]{Tor38} Tornheim, L.: Sums of $n$-th powers in fields of prime characteristic, {\it Duke Math. J.} {\bf 4} (1938), no. 2, 359--362.

\bibitem[TZS19]{TZS19} Tang, G.; Zhou, Y.; Su, H.: Matrices over a commutative ring as sums of three idempotents or three involutions, {\it Linear Multilinear Algebra} {\bf 67} (2019), no. 2, 267--277.

\bibitem[Vas03]{Vas03} Vasiu, A.: Surjectivity criteria for $p$-adic
representations. I, {\it Manuscripta Math.} {\bf 112} (2003), no. 3,
325--355.

\bibitem[Win98]{Win98} Winterhof, A.: On Waring's problem in finite fields, {\it Acta Arith.} {\bf 87} (1998), no. 2, 171--177.

\bibitem[Yo98]{Yo98} You, H.: Expressing a matrix over a finite field as the sum of two squares of invertible matrices, {\it J. Math. (Wuhan)} {\bf 18} (1998), no. 4, 383--388.

\end{thebibliography}
\end{document}